\documentclass[a4paper,reqno,12pt]{amsart}
\usepackage{amsmath,amsrefs,geometry}

\numberwithin{equation}{section}

\DeclareMathOperator{\p}{\overrightarrow{p}}
\DeclareMathOperator{\divergence}{div}
\DeclareMathOperator{\meas}{meas}
\DeclareMathOperator{\loc}{loc}
\DeclareMathOperator{\supp}{supp}
\DeclareMathOperator{\sgn}{sgn}
\renewcommand{\q}{\overrightarrow{q}}
\newcommand{\OO}{\text{O}}
\newcommand{\R}{\mathbb{R}}
\newcommand{\N}{\mathbb{N}}

\renewcommand{\[}{\left[}
\renewcommand{\]}{\right]}
\renewcommand{\(}{\left(}
\renewcommand{\)}{\right)}

\newtheorem{theorem}{Theorem}[section]
\newtheorem{lemma}[theorem]{Lemma}
\newtheorem{proposition}[theorem]{Proposition}
\newtheorem{remark}[theorem]{Remark}

\begin{document}

\date{October 2, 2014. {\it Final version:} July 22, 2015.}

\thanks{To appear in {\it Advances in Mathematics}}

\title[Critical anisotropic equations]{Decay estimates and a vanishing phenomenon for the solutions of critical anisotropic equations}

\author{J\'er\^ome V\'etois}

\address{J\'er\^ome V\'etois, Univ. Nice Sophia Antipolis, CNRS,  LJAD, UMR 7351, 06108 Nice, France.}
\email{vetois@unice.fr}

\begin{abstract}
We investigate the asymptotic behavior of solutions of anisotropic equations of the form $-\sum_{i=1}^n\partial_{x_i}\big(\left|\partial_{x_i}u\right|^{p_i-2}\partial_{x_i}u\big)=f\(x,u\)$ in $\R^n$, where $p_i>1$ for all $i=1,\dotsc,n$ and $f$ is a Caratheodory function with critical Sobolev growth. This problem arises in particular from the study of extremal functions for a class of anisotropic Sobolev inequalities. We establish decay estimates for the solutions and their derivatives, and we bring to light a vanishing phenomenon which occurs when the maximum value of the $p_i$ exceeds a critical value. 
\end{abstract}

\maketitle

\section{Introduction and main results}\label{Sec1}

We let $n\ge2$ and $\p=\(p_1,\dotsc,p_n\)$ be such that $p_i>1$ for all $i=1,\dotsc,n$ and $\sum_{i=1}^n1/p_i>1$. In this paper, we are interested in the solutions of problems of the form
\begin{equation}\label{Eq1}
\left\{\begin{aligned}&-\Delta_{\p} u=f\(x,u\)\quad\text{in }\R^n,\\
&u\in D^{1,\p}\(\R^n\),
\end{aligned}\right.
\end{equation}
where $\Delta_{\p} u:=\sum_{i=1}^n\partial_{x_i}\(\left|\partial_{x_i}u\right|^{p_i-2}\partial_{x_i}u\)$ is the anisotropic Laplace operator, $D^{1,\p}\(\R^n\)$ is the completion of $C^\infty_c\(\R^n\)$ with respect to the norm $\left\|u\right\|_{D^{1,\p}\(\R^n\)}:=\sum_{i=1}^n\(\int_{\R^n}\left|\partial_{x_i}u\right|^{p_i}dx\)^{1/p_i}$, and $f:\R^n\times\R\to\R$ is a Caratheodory function such that
\begin{equation}\label{Eq2}
\left|f\(x,s\)\right|\le\Lambda\left|s\right|^{p^*-1}\quad\text{for all }s\in\R\text{ and a.e. }x\in\R^n,
\end{equation}
for some real number $\Lambda>0$. Here, $p^*$ denotes the critical Sobolev exponent and is defined as
$$p^*:=\frac{n}{\sum_{i=1}^n\frac{1}{p_i}-1}=\frac{np}{n-p}\qquad\text{with}\quad\frac{1}{p}:=\frac{1}{n}\sum_{i=1}^n\frac{1}{p_i}\,.$$

\smallskip
The problem \eqref{Eq1} with $f\(x,u\)=\left|u\right|^{p^*-2}u$ appears in the study of extremal functions for a class of anisotropic Sobolev inequalities. Early references on anisotropic Sobolev inequalities are Nikol{\cprime}ski{\u\i}~\cite{Nik}, Troisi~\cite{Tro}, and Trudinger~\cite{Tru}. We also refer to Cianchi~\cite{Cia1} for a more recent work on the topic. Here we are interested in an inequality which appeared first in Troisi~\cite{Tro}.
%Anisotropic Sobolev inequalities were first studied by Nikol{\cprime}ski{\u\i}~\cite{Nik}. Here we are interested in an inequality which appeared in Troisi~\cite{Tro}. We also refer to Trudinger~\cite{Tru} for a more general class of inequalities, which in turn have been generalized by Cianchi~\cite{Cia1}. 
Among different equivalent versions (see Theorem~\ref{Th4} below), this inequality can be stated as
\begin{equation}\label{Eq3}
\int_{\R^n}\left|u\right|^{p^*}dx\le C\(\sum_{i=1}^n\int_{\R^n}\left|\partial_{x_i}u\right|^{p_i}dx\)^{p^*/p}
\end{equation}
for some constant $C=C\(n,\p\)$ and for all functions $u\in C^\infty_c\(\R^n\)$. The inequality \eqref{Eq3} enjoys an anisotropic scaling law (see \eqref{Eq9} below). As a corollary of the work of El~Hamidi--Rakotoson~\cite{ElHRak}, we obtain in Theorem~\ref{Th5} below that there exist extremal functions for the inequality \eqref{Eq3} provided that $p_i<p^*$ for all $i=1,\dotsc,n$. 

\smallskip
In the presence of anisotropy, namely when the $p_i$ are not all equal, there is no explicit formula for the extremal functions of \eqref{Eq3}. This motivates to find a priori estimates for these functions, and more generally for the solutions of equations of type \eqref{Eq1}. The main difficulties in this work come from the non-homogeneity of the problem and the lack of radial symmetry. 

\medskip
As a more general motivation, the solutions of problems of type \eqref{Eq1} with $f\(x,u\)=\left|u\right|^{p^*-2}u$ turn out to play a central role in the blow-up theories of critical equations in general domains. Possible references in book form on this subject and its applications in the isotropic regime are Druet--Hebey--Robert~\cite{DruHebRob}, Ghoussoub~\cite{Gho}, and Struwe~\cite{Str2}. A first step in the direction of a blow-up theory in the anisotropic regime was taken in El~Hamidi--V\'etois~\cite{ElHVet} where we extended the bubble tree decompositions of Struwe~\cite{Str1}. Now, if one wants to go further and investigate a pointwise blow-up theory, then it is essential to know the asymptotic behavior of the solutions of \eqref{Eq1} with $f\(x,u\)=\left|u\right|^{p^*-2}u$. The results in this paper can be seen as a crucial step in this direction.

\medskip
Anisotropic equations of type \eqref{Eq1} have received much attention in recent years. In addition to the above cited references~\cites{ElHRak,ElHVet} and without pretending to be exhaustive, we mention for instance the works by Cianchi~\cite{Cia2} on symmetrization properties, C\^irstea--V\'etois~\cite{CirVet} on the fundamental solutions, Cupini--Marcellini--Mascolo~\cite{CupMarMas} on the local boundedness of solutions, Fragal\`a--Gazzola--Kawohl~\cite{FraGazKaw} on the existence and non-existence of solutions in bounded domains, Lieberman~\cite{Lie} on gradient estimates, Namlyeyeva--Shishkov--Skrypnik~\cite{NamShiSkr} on singular solutions, and V\'etois~\cite{Vet2} on vanishing properties of solutions. More references can be found for instance in~\cites{Vet2}.

\medskip
Throughout this paper, we denote 
\begin{equation}\label{Eq4}
p_+:=\max\(\left\{p_i\in\p\right\}\)\quad\text{and}\quad p_*:=\frac{n-1}{\sum_{i=1}^n\frac{1}{p_i}-1}=\frac{p\(n-1\)}{n-p}\,.
\end{equation}
The exponent $p_*$ is known to play a critical role in several results on the asymptotic behavior of solutions of second order elliptic equations (see the historic paper of Serrin~\cite{Ser}, see also for instance the more recent paper of Serrin--Zou~\cite{SerZou} and the references therein). 

\medskip
Our first result is as follows.

\begin{theorem}\label{Th1}
Assume that $p_+<p_*$. Let $f:\R^n\times\R\to\R$ be a Caratheodory function such that \eqref{Eq2} holds true and $u$ be a solution of \eqref{Eq1}. Then there exists a constant $C_0=C_0\(n,\p,\Lambda,u\)$ such that 
\begin{equation}\label{Th1Eq1}
\left|u\(x\)\right|^{p_*}+\sum_{i=1}^n\left|\partial_{x_i}u\(x\)\right|^{p_i}\le C_0\bigg(1+\sum_{i=1}^n\left|x_i\right|^{\frac{p_*p_i}{p_*-p_i}}\bigg)^{-1}\quad\text{for a.e. }x\in\R^n,
\end{equation} 
where $p_*$ is as in \eqref{Eq4}.
\end{theorem}

We point out that the decay rate in \eqref{Th1Eq1} is the same as the one obtained in C\^irstea--V\'etois~\cite{CirVet} for the fundamental solutions in $\R^n$, namely the solutions of the equation $-\Delta_{\p}u=\delta_0$ in $\R^n$, where $\delta_0$ is the Dirac mass at the point 0.

\medskip
In case all $p_i$ are equal to $p$, as part of a more general result, Alvino--Ferone--Trombetti--Lions~\cite{AlvFerTroLio} proved that the best constant in the inequality \eqref{Eq3} is attained by the functions
\begin{equation}\label{Eq5}
u_{a,b}\(x\):=\bigg(a+b\sum_{i=1}^n\left|x_i\right|^{\frac{p}{p-1}}\bigg)^{\frac{p-n}{p}}
\end{equation}
for all $a,b>0$. Moreover, Cordero-Erausquin--Nazaret--Villani~\cite{CorNazVil} proved that the functions \eqref{Eq5} are the only extremal functions of \eqref{Eq3}. In case where the norm of the gradient in \eqref{Eq3} is replaced bv the Euclidean norm, the existence of radially symmetric extremal functions was found by Aubin~\cite{Aub}, Rodemich~\cite{Rod}, and Talenti~\cite{Tal}. Since $p/\(p_*-p\)=\(n-p\)/\(p-1\)$, the decay rate in \eqref{Eq5} coincides with the one in \eqref{Th1Eq1}.

\medskip
In case of the Laplace operator ($p_i=2$), Caffarelli--Gidas--Spruck~\cite{CafGidSpr} (see also Chen--Li~\cite{ChenLi}) proved that every positive solution of \eqref{Eq1} with $f\(x,u\)=u^{p^*-1}$ is of the form \eqref{Eq5}. This result can be extended to the case where all $p_i$ are equal to $p\in\(1,n\)$ for positive solutions satisfying the one-dimensional symmetry $u\(x\)=u\big(\sum_{i=1}^n\left|x_i\right|^{\frac{p}{p-1}}\big)$ for all $x\in\R^n$. Indeed, this result has been proved by  Guedda--V\'eron~\cite{GueVer} in case of positive, radially symmetric solutions for the $p$--Laplace equation $-\divergence\(\left|\nabla u\right|^{p-2}\nabla u\)=u^{p^*-1}\quad\text{in }\R^n$, and it can easily be seen that both cases lead to the same ordinary differential equation. We also mention that radial symmetry results have been established for positive solutions in $D^{1,p}\(\R^n\)$ in the case of $p$--Laplace equations (see Damascelli--Merch\'an--Montoro--Sciunzi~\cite{DamMerMonSci}, Damascelli--Ramaswamy~\cite{DamRam}, Sciunzi~\cite{Sci}, and V\'etois~\cite{Vet3}). 

\medskip
Theorem~\ref{Th1} has been proved in V\'etois~\cite{Vet3} in case of the $p$--Laplace operator. We also refer in case of the Laplace operator ($p_i=2$) to Jannelli--Solimini~\cite{JanSol}, where the decay estimate \eqref{Th1Eq1} has been proved to hold true for solutions of \eqref{Eq1}  with right-hand side $f\(x,u\)=\sum_{i=1}^Na_i\(x\)\left|u\right|^{q_i^*-2}u$, where $q_i^*:=2^*\(1-1/q_i\)$, $q_i\in\(n/2,\infty\]$, $\left|a_i\(x\)\right|=\OO\big(\left|x\right|^{-n/q_i}\big)$ for large $\left|x\right|$, and $a_i$ belongs to the Marcinkiewicz space $M^{q_i}\(\R^n\)$ for all $i=1,\dotsc,N$.

\medskip
The next results concern the case $p_+\ge p_*$, namely $p_i\ge p_*$ for some index $i$. In particular, we are now exclusively in the case where the exponents $p_i$ are not all equal.

\medskip
In the limit case $p_+=p_*$, we prove the following result.

\begin{theorem}\label{Th2}
Assume that $p_+=p_*$. Let $f:\R^n\times\R\to\R$ be a Caratheodory function such that \eqref{Eq2} holds true and $u$ be a solution of \eqref{Eq1}. Then for any $q>p_*$, there exists a constant $C_q=C\(n,\p,\Lambda,u,q\)$ such that 
\begin{equation}\label{Th2Eq}
\left|u\(x\)\right|^{q}+\sum_{i=1}^n\left|\partial_{x_i}u\(x\)\right|^{p_i}\le C_q\bigg(1+\sum_{i=1}^n\left|x_i\right|^{\frac{qp_i}{q-p_i}}\bigg)^{-1}\quad\text{for a.e. }x\in\R^n.
\end{equation}
\end{theorem}

Beyond this limit case, namely when $p_*<p_+<p^*$, we find the following result.

\begin{theorem}\label{Th3}
Assume that $p_*<p_+<p^*$. Let $f:\R^n\times\R\to\R$ be a Caratheodory function such that \eqref{Eq2} holds true and $u$ be a solution of \eqref{Eq1}. Then there exist a real number $q_0=q_0\(n,\p\)<p_+$ such that the two following assertions hold true.
\renewcommand{\labelenumi}{(\roman{enumi})}
\begin{enumerate}
\item There exists a constant $R_0=R_0\(n,\p,\Lambda,u\)$ such that 
\begin{equation}\label{Th3Eq1}
u\(x\)=0\quad\text{for all }x\in\R^n\text{ such that }\sum_{i\in\mathcal{I}_0}\left|x_i\right|\ge R_0\,,
\end{equation}
where $\mathcal{I}_0$ is the set of all indices $i$ such that $p_i>q_0$. Moreover, $\mathcal{I}_0\ne\emptyset$ due to $q_0<p_+$. 
\item For any $q>q_0$, there exists a constant $C_q=C\(n,\p,\Lambda,u,q\)$ such that
\begin{equation}\label{Th3Eq2}
\left|u\(x\)\right|^{q}+\sum_{i=1}^n\left|\partial_{x_i}u\(x\)\right|^{p_i}\le C_q\bigg(1+\sum_{i\in\mathcal{I}_0^c}\left|x_i\right|^{\frac{qp_i}{q-p_i}}\bigg)^{-1}\quad\text{for a.e. }x\in\R^n,
\end{equation}
where $\mathcal{I}_0^c:=\left\{1,\dotsc,n\right\}\backslash\mathcal{I}_0$.
\end{enumerate}
\end{theorem}

We are able, moreover, to give an explicit definition in terms of $n$ and $\p$ of a real number $q_0$ satisfying the above result (see Section~\ref{Sec7}).

\medskip
The dependence on $u$ of the constants $C_0$, $C_q$, and $R_0$ in the above results will be made more precise in Remarks~\ref{Rem2} and~\ref{Rem3}.

\medskip
As a remark about the support of solutions, by a result in V\'etois~\cite{Vet2}, we have that for any nonnegative solution $u$ of \eqref{Eq1} with $f\(x,u\)$ as in \eqref{Eq2} (see~\cite{Vet2} for the general assumptions), if $u\(x\)=0$ for some $x\in\R^n$, then we have $u\equiv0$ on the affine subspace $\left\{y\in\R^n:\,y_i=x_i\quad\forall i=\left\{1,\dotsc,n\right\}\backslash\mathcal{I}_-\right\}$, where $\mathcal{I}_-$ is the set of all indices $i$ such that $p_i=\min\(\left\{p_j\in\p\right\}\)$. 
In case all $p_i$ are equal to $p$, we obtain that either $u>0$ or $u\equiv0$, thus recovering the same result as Vazquez~\cite{Vaz} found for the $p$--Laplace operator. In the presence of anisotropy, as shows for instance Theorem~\ref{Th3}, this result does not hold true in general on the whole $\R^n$. 

\medskip
We also point out that in the limit case $p_+=p^*$,  we are able to construct quasi-explicit examples of solutions of \eqref{Eq1} with $f\(x,u\)=\left|u\right|^{p^*-2}u$ for anisotropic configurations of type $\p=\(p_-,\dotsc,p_-,p_+,\dotsc,p_+\)$ by using the method of separation of variables (see V\'etois~\cite{Vet1}). These solutions turn out to vanish in the $i$-th directions corresponding to $p_i=p_+$, exactly like what we prove to be true in Theorem~\ref{Th3} in case $p_*<p_+<p^*$.

\medskip
The paper is organized as follows. In Section~\ref{Sec2}, we present different equivalent versions of the anisotropic Sobolev inequality, and we study the existence and scaling properties of extremal functions for these inequalities.

\medskip
Section~\ref{Sec3} is concerned with preliminary properties satisfied by the solutions of \eqref{Eq1}, namely global boundedness results and a weak decay estimate.

\medskip
In Sections \ref{Sec4} and \ref{Sec5}, we perform a Moser-type iteration scheme inspired from the one developed in C\^irstea--V\'etois~\cite{CirVet} for the fundamental solutions. In order to treat a large part of the proofs in a unified way, we consider a general family of domains defined as
\begin{equation}\label{Eq6}
\Omega_{\q}\(\mathcal{I}_1,R_1,\mathcal{I}_2,R_2,\lambda\):=\bigg\{x\in\R^n:\sum_{i\in\mathcal{I}_1}\left|x_i\right|^{q_i}\hspace{-2pt}<\hspace{-2pt}\(1+\lambda\)R_1\text{ and }\bigg|\sum_{i\in\mathcal{I}_2}\left|x_i\right|^{q_i}-R_2\bigg|\hspace{-2pt}<\hspace{-2pt}\lambda R_2\bigg\},
\end{equation}
where $\lambda\in\(0,1\)$, $R_1,R_2>0$, $\mathcal{I}_1$ and $\mathcal{I}_2$ are two disjoint subsets of $\left\{1,\dotsc,n\right\}$, $\mathcal{I}_2\ne\emptyset$, and $\q=\(q_i\)_{i\in\mathcal{I}_1\cup\mathcal{I}_2}$ is such that $q_i>1$ for all $i\in\mathcal{I}_1\cup\mathcal{I}_2$. On these domains, we prove that the solutions of \eqref{Eq1} satisfy reverse H\"older-type inequalities of the form
\begin{equation}\label{Eq7}
\left\|u\right\|_{L^\gamma\(\Omega_{\q}\(\mathcal{I}_1,R_1,\mathcal{I}_2,R_2,\lambda\)\)}^\gamma\le C\max_{i\in\mathcal{I}_1\cup\mathcal{I}_2}\Big(\(\lambda'-\lambda\)^{-p_i}\\
R_{\delta_i}^{-\frac{p_i}{q_i}}\left\|u\right\|_{L^{\gamma_i}\(\Omega_{\q}\(\mathcal{I}_1,R_1,\mathcal{I}_2,R_2,\lambda'\)\)}^{\gamma_i}\Big)^\frac{n}{n-p}
\end{equation}
for all $\gamma>p_*-1$ and $\lambda<\lambda'\in\(0,1/2\]$, where $\delta_i:=1$ if $i\in\mathcal{I}_1$, $\delta_i:=2$ if $i\in\mathcal{I}_2$, $\gamma_i:=\frac{n-p}{n}\gamma+p_i-p$ for all $i\in\mathcal{I}_1\cup\mathcal{I}_2$,  and $C=C\(n,\p,\q,u,\gamma\)$ (see Lemma~\ref{Lem4} for more details on the dependence of the constant with respect to  $u$ and $\gamma$). Since the right-hand side of \eqref{Eq7} involves different exponents $\gamma_i$ in the anisotropic case, the number of exponents in the estimates may grow exponentially when iterating this inequality. We overcome this issue in Section~\ref{Sec5} by controlling the values of the exponents with respect to the number of iterations. 

\medskip
In Section~\ref{Sec6}, we prove a vanishing result which will give Point (i) in Theorem~\ref{Th3}. We prove this result by applying our iteration scheme with $R_1=R_2^{1/\varepsilon}$ for small real numbers $\varepsilon>0$ and $\mathcal{I}_1$, $\mathcal{I}_2$ being the sets of all indices $i$ such that $p_i<p_0$, $p_i=p_0$, respectively, for some large enough real number $p_0\in\p$ (see \eqref{Lem8Eq1} for the exact condition on $p_0$). Passing to the limit into our iteration scheme, we obtain a pointwise estimate of the form
\begin{equation}\label{Eq8}
\left\|u\right\|_{L^\infty(\Omega_{\q}(\mathcal{I}_1,R^{1/\varepsilon},\mathcal{I}_2,R,1/4))}\le\big(CR^{-\frac{1}{p_\varepsilon}}\big)^{\frac{1}{\varepsilon}}
\end{equation}
for some constant $C=C\(n,\p,u\)$ (see Lemma~\ref{Lem8}). When $R$ is large enough, the right-hand side of \eqref{Eq8} converges to 0 as $\varepsilon\to0$, and we thus obtain our vanishing result. 

\medskip
In Section~\ref{Sec7}, we prove Theorems~\ref{Th1},~\ref{Th2}, and we complete the proof of Theorem~\ref{Th3} by proving the decay estimates \eqref{Th3Eq2}. The proofs of these results rely again on our iteration scheme, this time applied with $\mathcal{I}_1=\emptyset$ and $\mathcal{I}_2$ being the set of all indices $i$ such that $p_i\le\overline{p}_0$ for some real number $\overline{p}_0$ (see \eqref{Eq15}). 

\medskip
Finally, in Appendix~\ref{App}, we prove a weak version of Kato's inequality which is used in Sections~\ref{Sec3} and~\ref{Sec4}.

\medskip\noindent
{\bf Acknowledgments.} The author wishes to express his gratitude to Fr\'ed\'eric Robert for helpful comments on the manuscript.

\section{Application to the extremal functions of a class of\\anisotropic Sobolev inequalities}\label{Sec2}

As mentioned in the introduction, one of our main motivation in this paper is to apply our results to the extremal functions of a class of anisotropic Sobolev inequalities which originates from Troisi~\cite{Tro}. In this section, we first present in Theorem~\ref{Th4} below different equivalent versions of these inequalities, and we then prove in Theorem~\ref{Th5} that all these inequalities have extremal functions, and that with a suitable change of scale, these extremal functions are solutions of \eqref{Eq1} with $f\(x,u\)=\left|u\right|^{p^*-2}u$.

\medskip
We state the equivalent versions of the anisotropic Sobolev inequalities as follows.

\begin{theorem}\label{Th4}
The following inequalities hold true.
\renewcommand{\labelenumi}{(\roman{enumi})}
\begin{enumerate}
\item There exists a constant $C=C\big(n,\p\big)$ such that
\begin{equation}\label{Th4Eq1}
\(\int_{\R^n}\left|u\right|^{p^*}dx\)^{n/p^*}\le C\prod_{i=1}^n\(\int_{\R^n}\left|\partial_{x_i}u\right|^{p_i}dx\)^{1/p_i}\quad\forall u\in C^\infty_c\(\R^n\).
\end{equation}
\item  For any $\overrightarrow{\theta}=\(\theta_1,\dotsc,\theta_n\)$ such that $\theta_i>0$ for all $i=1,\dotsc,n$ and $\sum_{i=1}^n1/\theta_i=n/p$, there exists a constant $C_{\overrightarrow{\theta}}=C\big(n,\p,\overrightarrow{\theta}\big)$ such that
\begin{equation}\label{Th4Eq2}
\(\int_{\R^n}\left|u\right|^{p^*}dx\)^{p/p^*}\le C_{\overrightarrow{\theta}}\sum_{i=1}^n\(\int_{\R^n}\left|\partial_{x_i}u\right|^{p_i}dx\)^{\theta_i/p_i}\quad\forall u\in C^\infty_c\(\R^n\).
\end{equation}
In particular, we get \eqref{Eq3} in case $\theta_i=p_i$ for all $i=1,\dotsc,n$.
\end{enumerate}
\end{theorem}

As a remark, the inequalities \eqref{Th4Eq1} and \eqref{Th4Eq2} enjoy an anisotropic scaling law. Indeed, it can easily be seen that every integral in these inequalities are invariant with respect to the change of scale $u\mapsto u_\lambda$, where 
\begin{equation}\label{Eq9}
u_\lambda\(x\)=\lambda u\(\lambda^{\(p^*-p_1\)/p_1}x_1,\dotsc,\lambda^ {\(p^*-p_n\)/p_n}x_n\)
\end{equation}
for all $\lambda>0$ and $x\in\R^n$. 

\proof[Proof of Theorem~\ref{Th4}]
We refer to Troisi~\cite{Tro}*{Theorem~1.2} for the proof of the inequality \eqref{Th4Eq1}. Then the inequality \eqref{Th4Eq2} follows from \eqref{Th4Eq1} by applying an inequality of weighted arithmetic and geometric means. As a remark, we can also obtain \eqref{Th4Eq1} from \eqref{Th4Eq2} by applying the change of scale \eqref{Pr1Eq2} below.
\endproof

Regarding the extremal functions of \eqref{Th4Eq1} and \eqref{Th4Eq2}, we prove the following result. The existence part in this result will be obtained as a corollary of the work of El~Hamidi--Rakotoson~\cite{ElHRak} and Proposition~\ref{Pr1} below.

\begin{theorem}\label{Th5}
If $p_+<p^*$, then there exist extremal functions $u\in D^{1,\p}\(\R^n\)$, $u\ne0$, of  \eqref{Th4Eq1} and \eqref{Th4Eq2}. Moreover, for any extremal function $u$ of \eqref{Th4Eq1} or \eqref{Th4Eq2}, there exist $\mu_1,\dotsc,\mu_n>0$ such that the function $x\in\R^n\mapsto u\(\mu_1x_1,\dotsc,\mu_nx_n\)$ is a constant-sign solution of \eqref{Eq1} with $f\(x,u\)=\left|u\right|^{p^*-2}u$. In particular, every extremal function of \eqref{Th4Eq1} or \eqref{Th4Eq2} satisfies the a priori estimates in Theorems~\ref{Th1},~\ref{Th2}, and~\ref{Th3}.
\end{theorem}

As a remark, due to the scaling law \eqref{Eq9}, every extremal function of \eqref{Th4Eq1} or \eqref{Th4Eq2} generates in fact an infinite family of extremal functions.

\medskip
Preliminary to the proof of Theorem~\ref{Th5}, we prove the following result.

\begin{proposition}\label{Pr1}
Let $\overrightarrow{\theta}=\(\theta_1,\dotsc,\theta_n\)$ be such that $\theta_i>0$ for all $i=1,\dotsc,n$ and $\sum_{i=1}^n1/\theta_i=n/p$. Then the following assertions hold true.
\renewcommand{\labelenumi}{(\roman{enumi})}
\begin{enumerate}
\item For any extremal function $u$ of \eqref{Th4Eq2}, $u\circ\tau_{\overrightarrow{\theta}}$ is an extremal function of \eqref{Th4Eq1}, where
\begin{equation}\label{Pr1Eq1}
\tau_{\overrightarrow{\theta}}\(x\):=\left(\lambda_{\overrightarrow{\theta},1}x_1,\dotsc,\lambda_{\overrightarrow{\theta},n}x_n\right),\quad\lambda_{\overrightarrow{\theta},i}:=\theta_i^{1/\theta_i}\prod_{j=1}^ n\theta_j^{-p/\(n\theta_i\theta_j\)}\,,
\end{equation}
for all $x\in\R^n$ and $i=1,\dotsc,n$.
\item For any extremal function $u$ of \eqref{Th4Eq1}, $u\circ\sigma_{\overrightarrow{\theta},u}\circ\tau_{\overrightarrow{\theta}}^{-1}$ is an extremal function of \eqref{Th4Eq2}, where $\tau_{\overrightarrow{\theta}}$ is as in \eqref{Pr1Eq1} and
\begin{equation}\label{Pr1Eq2}
\sigma_{\overrightarrow{\theta},u}\(x\):=\left(\mu_{\overrightarrow{\theta},1}\(u\)x_1,\dotsc,\mu_{\overrightarrow{\theta},n}\(u\)x_n\right),\quad\mu_{\overrightarrow{\theta},i}\(u\):=\frac{\overset{n}{\underset{j=1}{\prod}}\(\int_{\R^n}\left|\partial_{x_j}u\right|^{p_j}dx\)^{p/\(n\theta_ip_j\)}}{\(\int_{\R^n}\left|\partial_{x_i}u\right|^{p_i}dx\)^{1/p_i}}\,,
\end{equation}
for all $x\in\R^n$ and $i=1,\dotsc,n$.
\end{enumerate}
\end{proposition}

\proof[Proof of Proposition~\ref{Pr1}]
We begin with proving Point (i). We fix an extremal function $u_0$ of \eqref{Th4Eq2}. Since $\sum_{i=1}^n1/\theta_i=n/p$, we obtain
\begin{equation}\label{Pr1Eq3}
\prod_{i=1}^n\(\int_{\R^n}\big|\partial_{x_i}\big(u_0\circ\tau_{\overrightarrow{\theta}}\big)\big|^{p_i}dx\)^{p/\(np_i\)}\le\frac{p}{n}\sum_{i=1}^n\frac{1}{\theta_i}\(\int_{\R^n}\big|\partial_{x_i}\big(u_0\circ\tau_{\overrightarrow{\theta}}\big)\big|^{p_i}dx\)^{\theta_i/p_i}.
\end{equation}
For any function $u\in D^{1,\p}\(\R^n\)$, simple calculations give
\begin{equation}\label{Pr1Eq4}
\sum_{i=1}^n\frac{1}{\theta_i}\(\int_{\R^n}\big|\partial_{x_i}\big(u\circ\tau_{\overrightarrow{\theta}}\big)\big|^{p_i}dx\)^{\theta_i/p_i}=\(\prod_{j=1}^ n\theta_j^{-p/\(n\theta_j\)}\)\sum_{i=1}^n\(\int_{\R^n}\left|\partial_{x_i}u\right|^{p_i}dx\)^{\theta_i/p_i}
\end{equation}
and
\begin{equation}\label{Pr1Eq5}
\int_{\R^n}\big|u\circ\tau_{\overrightarrow{\theta}}\big|^{p^*}dx=\int_{\R^n}\left|u\right|^{p^*}dx\,.
\end{equation}
By invertibility of $\tau_{\overrightarrow{\theta}}$ and since $u_0$ is an extremal function of \eqref{Th4Eq2}, it follows from \eqref{Pr1Eq4} and \eqref{Pr1Eq5} that
\begin{equation}\label{Pr1Eq6}
\frac{\overset{n}{\underset{i=1}{\sum}}\frac{1}{\theta_i}\(\int_{\R^n}\big|\partial_{x_i}\big(u_0\circ\tau_{\overrightarrow{\theta}}\big)\big|^{p_i}dx\)^{\theta_i/p_i}}{\big(\int_{\R^n}\big|u_0\circ\tau_{\overrightarrow{\theta}}\big|^{p^*}dx\big)^{p/p^*}}=\inf_{\underset{u\ne0}{u\in D^{1,\p}\(\R^n\)}}\frac{\overset{n}{\underset{i=1}{\sum}}\frac{1}{\theta_i}\(\int_{\R^n}\left|\partial_{x_i}u\right|^{p_i}dx\)^{\theta_i/p_i}}{\left(\int_{\R^n}\left|u\right|^{p^*}dx\right)^{p/p^*}}\,.
\end{equation}
Now, we claim that
\begin{equation}\label{Pr1Eq7}
\inf_{\underset{u\ne0}{u\in D^{1,\p}\(\R^n\)}}\frac{\overset{n}{\underset{i=1}{\sum}}\frac{1}{\theta_i}\(\int_{\R^n}\left|\partial_{x_i}u\right|^{p_i}dx\)^{\theta_i/p_i}}{\left(\int_{\R^n}\left|u\right|^{p^*}dx\right)^{p/p^*}}\le\frac{n}{p}\cdot\inf_{\underset{u\ne0}{u\in D^{1,\p}\(\R^n\)}}\frac{\overset{n}{\underset{i=1}{\prod}}\(\int_{\R^n}\left|\partial_{x_i}u\right|^{p_i}dx\)^{p/\(np_i\)}}{\left(\int_{\R^n}\left|u\right|^{p^*}dx\right)^{p/p^*}}\,.
\end{equation}
We prove this claim. For any function $u\in D^{1,\p}\(\R^n\)$, $u\ne0$, by applying the change of scale \eqref{Pr1Eq2}, we obtain
\begin{equation}\label{Pr1Eq8}
\(\int_{\R^n}\big|\partial_{x_i}\big(u\circ\sigma_{\overrightarrow{\theta},u}\big)\big|^{p_i}dx\)^{\theta_i/p_i}=\prod_{j=1}^n\(\int_{\R^n}\left|\partial_{x_j}u\right|^{p_j}dx\)^{p/\(np_j\)}
\end{equation}
for all $i=1,\dotsc,n$, and
\begin{equation}\label{Pr1Eq9}
\int_{\R^n}\big|u\circ\sigma_{\overrightarrow{\theta},u}\big|^{p^*}dx=\int_{\R^n}\left|u\right|^{p^*}dx\,.
\end{equation}
Since $\sum_{i=1}^n1/\theta_i=n/p$, it follows from \eqref{Pr1Eq8} and \eqref{Pr1Eq9} that
\begin{equation}\label{Pr1Eq10}
\frac{\overset{n}{\underset{i=1}{\sum}}\frac{1}{\theta_i}\(\int_{\R^n}\big|\partial_{x_i}\big(u\circ\sigma_{\overrightarrow{\theta},u}\big)\big|^{p_i}dx\)^{\theta_i/p_i}}{\left(\int_{\R^n}\big|u\circ\sigma_{\overrightarrow{\theta},u}\big|^{p^*}dx\right)^{p/p^*}}=\frac{n}{p}\cdot\frac{\overset{n}{\underset{i=1}{\prod}}\(\int_{\R^n}\left|\partial_{x_i}u\right|^{p_i}dx\)^{p/\(np_i\)}}{\left(\int_{\R^n}\left|u\right|^{p^*}dx\right)^{p/p^*}}\,,
\end{equation}
and hence we obtain \eqref{Pr1Eq7}. It follows from \eqref{Pr1Eq3}, \eqref{Pr1Eq6} and \eqref{Pr1Eq7} that $u_0\circ\tau_{\overrightarrow{\theta}}$ is an extremal function of \eqref{Th4Eq1}. This ends the proof of Point (i).

Now, we prove Point (ii). We fix an extremal function $u_0$ of \eqref{Th4Eq1}. By \eqref{Pr1Eq10} and since $\sum_{i=1}^n1/\theta_i=n/p$, we obtain
\begin{align}\label{Pr1Eq11}
\frac{\overset{n}{\underset{i=1}{\sum}}\frac{1}{\theta_i}\(\int_{\R^n}\big|\partial_{x_i}\big(u_0\circ\sigma_{\overrightarrow{\theta},u_0}\big)\big|^{p_i}dx\)^{\theta_i/p_i}}{\left(\int_{\R^n}\big|u_0\circ\sigma_{\overrightarrow{\theta},u_0}\big|^{p^*}dx\right)^{p/p^*}}&=\frac{n}{p}\cdot\inf_{\underset{u\ne0}{u\in D^{1,\p}\(\R^n\)}}\frac{\overset{n}{\underset{i=1}{\prod}}\(\int_{\R^n}\left|\partial_{x_i}u\right|^{p_i}dx\)^{p/\(np_i\)}}{\left(\int_{\R^n}\left|u\right|^{p^*}dx\right)^{p/p^*}}\nonumber\\
&\le\inf_{\underset{u\ne0}{u\in D^{1,\p}\(\R^n\)}}\frac{\overset{n}{\underset{i=1}{\sum}}\frac{1}{\theta_i}\(\int_{\R^n}\left|\partial_{x_i}u\right|^{p_i}dx\)^{\theta_i/p_i}}{\left(\int_{\R^n}\left|u\right|^{p^*}dx\right)^{p/p^*}}\,.
\end{align}
It follows from \eqref{Pr1Eq4} and \eqref{Pr1Eq11} that $u_0\circ\sigma_{\overrightarrow{\theta},u_0}\circ\tau_{\overrightarrow{\theta}}^{-1}$ is an extremal function of \eqref{Th4Eq2}. This ends the proof of Point (ii).
\endproof

Now, we can prove Theorem~\ref{Th5} by using Proposition~\ref{Pr1}.

\proof[Proof of Theorem~\ref{Th5}]
We prove the results for the sole inequality \eqref{Eq3}. The results for \eqref{Th4Eq1} and \eqref{Th4Eq2} then follow from Proposition~\ref{Pr1}.

First, in case $p_+<p^*$, the existence of extremal functions of \eqref{Eq3} follows from the work of El~Hamidi--Rakotoson~\cite{ElHRak}. Indeed, it has been proven in~\cite{ElHRak} that there exist minimizers for
\begin{equation}\label{Th5Eq1}
\mathcal{I}:=\inf_{\underset{\int_{\R^n}\left|u\right|^{p^*}dx=1}{u\in D^{1,\p}\(\R^n\)}}\sum_{i=1}^n\frac{1}{p_i}\int_{\R^n}\left|\partial_{x_i}u\right|^{p_i}dx\,.
\end{equation}
This infimum is connected with \eqref{Eq3} by the change of scale $u\mapsto\mu_{\p,u}^{-1}\cdot u\circ\rho_{\p,u}$, where
$$\rho_{\p,u}\(x\):=\mu_{\p,u}\cdot\tau_{\p}\(x\),\quad\mu_{\p,u}:=\(\int_{\R^n}\left|u\right|^{p^*}dx\)^{p/\(np^*\)},$$
and $\tau_{\p}\(x\)$ is as in \eqref{Pr1Eq1} for all $x\in\R^n$ and $u\in D^{1,\p}\(\R^n\)$, $u\ne0$. More precisely, simple calculations give
\begin{equation}\label{Th5Eq2}
\sum_{i=1}^n\frac{1}{p_i}\int_{\R^n}\big|\partial_{x_i}\big(\mu_{\p,u}^{-1}\cdot u\circ\rho_{\p,u}\big)\big|^{p_i}dx=\(\prod_{j=1}^np_j^{-p/\(np_j\)}\)\frac{\overset{n}{\underset{i=1}{\sum}}\int_{\R^n}\big|\partial_{x_i}u\big|^{p_i}dx}{\(\int_{\R^n}\left|u\right|^{p^*}dx\)^{p/p^*}}
\end{equation}
and
\begin{equation}\label{Th5Eq3}
\int_{\R^n}\big|\mu_{\p,u}^{-1}\cdot u\circ\rho_{\p,u}\big|^{p^*}dx=1\,,
\end{equation}
and hence
\begin{equation}\label{Th5Eq4}
\mathcal{I}\le\(\prod_{j=1}^np_j^{-p/\(np_j\)}\)\inf_{\underset{u\ne0}{u\in D^{1,\p}\(\R^n\)}}\frac{\overset{n}{\underset{i=1}{\sum}}\int_{\R^n}\big|\partial_{x_i}u\big|^{p_i}dx}{\(\int_{\R^n}\left|u\right|^{p^*}dx\)^{p/p^*}}\,.
\end{equation}
In particular, for any minimizer $u$ of \eqref{Th5Eq1}, since $\mu_{\p,u}=1$ and $\rho_{\p,u}=\tau_{\p}$, it follows from \eqref{Th5Eq2}--\eqref{Th5Eq4} that $u\circ\tau_{\p}^{-1}$ is an extremal function of \eqref{Eq3}. 

Next, we prove that the extremal functions of \eqref{Eq3} do not change sign. We let $C_0$ be the best constant and $u$ be an extremal function of \eqref{Eq3}. By writing $u=u_+-u_-$, where $u_+:=\max\(u,0\)$ and $u_-:=\max\(-u,0\)$, we obtain
\begin{multline}\label{Th5Eq5}
\sum_{i=1}^n\int_{\R^n}\left|\partial_{x_i}u\right|^{p_i}dx=\(\frac{1}{C_0}\int_{\R^n}\left|u\right|^{p^*}dx\)^{p/p_*}=\(\frac{1}{C_0}\int_{\R^n}u_-^{p^*}dx+\frac{1}{C_0}\int_{\R^n}u_+^{p^*}dx\)^{p/p_*}\\
\le\(\(\sum_{i=1}^n\int_{\R^n}\left|\partial_{x_i}u_-\right|^{p_i}dx\)^{p^*/p}+\(\sum_{i=1}^n\int_{\R^n}\left|\partial_{x_i}u_+\right|^{p_i}dx\)^{p^*/p}\)^{p/p_*}.
\end{multline}
It follows from \eqref{Th5Eq5} that either $u_-=0$ or $u_+=0$, and hence we obtain that the function $u$ has constant sign.

Finally, from the Euler-Lagrange equation satisfied by $u$, namely
$$-\sum_{i=1}^np_i\partial_{x_i}\(\left|\partial_{x_i}u\right|^{p_i-2}\partial_{x_i}u\)=\lambda\(u\)\left|u\right|^{p^*-2}u,\quad\text{where}\quad\lambda\(u\):=\frac{\overset{n}{\underset{i=1}{\sum}}p_i\int_{\R^n}\big|\partial_{x_i}u\big|^{p_i}dx}{\int_{\R^n}\left|u\right|^{p^*}dx}\,,$$
we derive that the function 
$x\in\R^n\mapsto\mu u\(\mu_1x_1,\dotsc,\mu_nx_n\)$ with $\mu_i:=\(\lambda\(u\)/p_i\)^{1/p_i}$ for all $i=1,\dotsc,n$ is a solution of \eqref{Eq1} with $f\(x,u\)=\left|u\right|^{p^*-2}u$. This ends the proof of Theorem~\ref{Th5}.
\endproof

\section{Preliminary results}\label{Sec3}

From now on, we are concerned with the general case of an arbitrary solution of \eqref{Eq1}. 

\medskip
For any $s\in\(0,\infty\)$ and any domain $\Omega\subset\R^n$, we define the weak Lebesgue space $L^{s,\infty}\(\Omega\)$ as the set of all measurable functions $u:\Omega\to\R$ such that
$$\left\|u\right\|_{L^{s,\infty}\(\Omega\)}:=\sup_{h>0}\big(h\cdot\meas\(\left\{\left|u\right|>h\right\}\)^{1/s}\big)<\infty\,,$$
where $\meas\(\left\{\left|u\right|>h\right\}\)$ is the measure of the set $\left\{x\in\Omega:\,\left|u\(x\)\right|>h\right\}$. The map $\left\|\cdot\right\|_{L^{s,\infty}\(\Omega\)}$ defines a quasi-norm on $L^{s,\infty}\(\Omega\)$. We refer, for instance, to the book of Grafakos~\cite{Gra} for the material on weak Lebesgue spaces. 

\medskip
The first result in this section is as follows.

\begin{lemma}\label{Lem1}
Assume that $p_+<p^*$. Let $f:\R^n\times\R\to\R$ be a Caratheodory function such that \eqref{Eq2} holds true. Then any solution of \eqref{Eq1} belongs to $W^{1,\infty}\(\R^n\)\cap L^{p_*-1,\infty}\(\R^n\)$, and hence by interpolation, to $L^s\(\R^n\)$ for all $s\in\(p_*-1,\infty\]$.
\end{lemma}

\proof[Proof of Lemma~\ref{Lem1}]
The $L^\infty$--boundedness of the solutions follows from a straightforward adaptation of El~Hamidi--Rakotoson~\cite{ElHRak}*{Propositions~1 and~2}, the first proposition being in turn adapted from Fragal\`a--Gazzola--Kawohl~\cite{FraGazKaw}*{Theorem~2}.

Once we have the $L^\infty$--boundedness of the solutions, we obtain the $L^\infty$--boundedness of the derivatives by applying Lieberman's gradient estimates~\cite{Lie}.

The proof of the $L^{p_*-1,\infty}$--boundedness of the solutions follows exactly the same arguments as in V\'etois~\cite{Vet3}*{Lemma~2.2}. One only has to replace $\left|\nabla u\right|^p$ by $\sum_{i=1}^n\left|\partial_{x_i}u\right|^{p_i}$.
\endproof

For any solution $u$ of \eqref{Eq2}, by Proposition~\ref{Pr2} in Appendix~\ref{App}, we obtain
\begin{equation}\label{Eq10}
-\Delta_{\p}\left|u\right|\le f\(x,u\)\cdot\sgn\(u\)\le\Lambda\left|u\right|^{p^*-1}\quad\text{in }\R^n,
\end{equation}
where $\sgn\(u\)$ denotes the sign of $u$ and the inequality is in the sense that for any nonnegative, smooth function $\varphi$ with compact support in $\R^n$, we have
$$\sum_{i=1}^n\int_{\R^n}\left|\partial_{x_i}\left|u\right|\right|^{p_i-2}\(\partial_{x_i}\left|u\right|\)\(\partial_{x_i}\varphi\)dx\le\Lambda\int_{\R^n}\left|u\right|^{p^*-1}\varphi\,dx\,.$$

\medskip
We prove the following result.

\begin{lemma}\label{Lem2}
For any real number $\Lambda>0$ and any nonnegative, nontrivial solution $v\in D^{1,p}\(\R^n\)$ of the inequality $-\Delta_{\p}v\le\Lambda v^{p^*-1}$ in $\R^n$, we have $\left\|v\right\|_{L^{p^*}\(\R^n\)}\ge\kappa_0$ for some constant $\kappa_0=\kappa_0\(n,p,\Lambda\)>0$.
\end{lemma}

\proof
By testing the inequality $-\Delta_{\p}v\le\Lambda v^{p^*-1}$ with the function $v$, and applying the anisotropic Sobolev inequality, we obtain
\begin{equation}\label{Lem2Eq}
\Lambda\int_{\R^n}v^{p^*}dx\ge\sum_{i=1}^n\int_{\R^n}\left|\partial_{x_i}v\right|^{p_i}dx\ge K\(\int_{\R^n}v^{p^*}dx\)^{\frac{n-p}{n}}
\end{equation}
for some constant $K=K\(n,\p\)$. The result then follows from \eqref{Lem2Eq} with $\kappa_0:=\(K/\Lambda\)^{\frac{n-p}{p^2}}$.
\endproof

As a last result in this section, we prove the following decay estimate. This result is not sharp, but it turns out to be a crucial ingredient in what follows.

\begin{lemma}\label{Lem3}
Assume that $p_+<p^*$. Let $\kappa_0$ be as in Lemma~\ref{Lem2}, $f:\R^n\times\R\to\R$ be a Caratheodory function such that \eqref{Eq2} holds true, and $u$ be a solution of \eqref{Eq1}. For any $\kappa>0$, we define
\begin{equation}\label{Lem3Eq1}
r_\kappa\(u\):=\inf\big(\big\{r>0\,:\,\left\|u\right\|_{L^{p^*}\(\R^n\backslash B_{\p}\(0,r\)\)}<\kappa\big\}\big),
\end{equation}
where $B_{\p}\(0,r\)$ is the open ball of center 0 and radius $r$ with respect to the distance function $d_{\p}$ defined as
\begin{equation}\label{Lem3Eq2}
d_{\p}\(x,y\):=\sum_{i=1}^n\left|x_i-y_i\right|^{\frac{\delta p_i}{p^*-p_i}}\quad\text{with}\quad\delta:=\frac{p^*-p_+}{p_+},
\end{equation}
for all $x,y\in\R^n$. Then for any $\kappa\in\(0,\kappa_0\)$ and $r>r_\kappa\(u\)$, there exists a constant $K_0=K_0\big(n,\p,\Lambda,\kappa,r,r_\kappa\(u\),\left\|u\right\|_{L^{p^*}\(\R^n\)}\big)$ such that
\begin{equation}\label{Lem3Eq3}
\left|u\(x\)\right|\le K_0\bigg(\sum_{i=1}^n\left|x_i\right|^{\frac{p_i}{p^*-p_i}}\bigg)^{-1}\quad\text{for all }x\in\R^n\backslash B_{\p}\(0,r\).
\end{equation}
\end{lemma}

\proof[Proof of Lemma~\ref{Lem3}]
This proof is adapted from V\'etois~\cite{Vet3}*{Lemma~3.1}
We fix $\Lambda>0$, $\kappa\in\(0,\kappa_0\)$, $\kappa'>\kappa_0$, $r>0$, and $r'\in\(0,r\)$. We claim that in order to obtain Lemma~\ref{Lem3}, it is sufficient to prove that there exists a constant $K_1=K_1\(n,\p,\kappa,\kappa',r,r'\)$ such that for any solution $u$ of \eqref{Eq1} such that $r_\kappa\(u\)\le r'$ and $\left\|u\right\|_{L^{p^*}\(\R^n\)}\le\kappa'$, we have
\begin{equation}\label{Lem3Eq4}
d_{\p}\(x,B_{\p}\(0,r''\)\)\left|u\(x\)\right|^\delta\le K_1\quad\text{for all }x\in\R^n\backslash B_{\p}\(0,r\),
\end{equation}
where $r'':=\(r+r'\)/2$. Indeed, for any  $x\in\R^n\backslash B_{\p}\(0,r\)$, we can write 
\begin{equation}\label{Lem3Eq5}
d_{\p}\(x,0\)\le d_{\p}\(x,B_{\p}\(0,r''\)\)+r''\le d_{\p}\(x,B_{\p}\(0,r''\)\)+\frac{r''}{r}d_{\p}\(x,0\),
\end{equation}
and hence by putting together \eqref{Lem3Eq4} and \eqref{Lem3Eq5}, we obtain
\begin{equation}\label{Lem3Eq6}
d_{\p}\(x,0\)\left|u\(x\)\right|^\delta\le\frac{r}{r-r''}\cdot K_1=\frac{2r}{r-r'}\cdot K_1.
\end{equation}
By definition of $d_{\p}$, \eqref{Lem3Eq3} then follows from \eqref{Lem3Eq6}. This proves our claim.

We prove \eqref{Lem3Eq4} by contradiction. Suppose that for any $\alpha\in\N$, there exists a Caratheodory function $f_\alpha:\R^n\times\R\to\R$ such that \eqref{Eq2} holds true, a solution $u_\alpha$ of \eqref{Eq1} with $f=f_\alpha$ such that $r_\kappa\(u_\alpha\)\le r'$ and $\left\|u_\alpha\right\|_{L^{p^*}\(\R^n\)}\le\kappa'$, and a point $x_\alpha\in\R^n\backslash B\(0,r\)$ such that
\begin{equation}\label{Lem3Eq7}
d_{\p}\(x_\alpha,B_{\p}\(0,r''\)\)\left|u_\alpha\(x_\alpha\)\right|^\delta>2\alpha\,.
\end{equation}
It follows from \eqref{Lem3Eq7} and Pol\'a\v{c}ik--Quittner--Souplet~\cite{PolQuiSou}*{Lemma~5.1} that there exists $y_\alpha\in\R^n\backslash B_{\p}\(0,r''\)$ such that
\begin{equation}\label{Lem3Eq8}
d_{\p}\(y_\alpha,B_{\p}\(0,r''\)\)\left|u_\alpha\(y_\alpha\)\right|^\delta>2\alpha\,,\quad\left|u_\alpha\(x_\alpha\)\right|\le\left|u_\alpha\(y_\alpha\)\right|,
\end{equation}
and 
\begin{equation}\label{Lem3Eq9}
\left|u_\alpha\(y\)\right|\le2^{1/\delta}\left|u_\alpha\(y_\alpha\)\right|\quad\text{for all }y\in B_{\p}\big(y_\alpha,\alpha\,\left|u_\alpha\(y_\alpha\)\right|^{-\delta}\big).
\end{equation}
For any $\alpha$ and $y\in\R^n$, we define
\begin{equation}\label{Lem3Eq10}
\widetilde{u}_\alpha\(y\):=\left|u_\alpha\(y_\alpha\)\right|^{-1}\cdot u_\alpha\(\tau_\alpha\(y\)\),
\end{equation}
where
$$\tau_\alpha\(y\):=y_\alpha+\big(\left|u_\alpha\(y_\alpha\)\right|^{\frac{p_1-p^*}{p_1}}y_1,\dots,\left|u_\alpha\(y_\alpha\)\right|^{\frac{p_n-p^*}{p_n}}y_n\big).$$
It follows from \eqref{Lem3Eq9} and \eqref{Lem3Eq10} that
\begin{equation}\label{Lem3Eq11}
\left|\widetilde{u}_\alpha\(0\)\right|=1\quad\text{and}\quad\left|\widetilde{u}_\alpha\(y\)\right|\le2^{1/\delta}\quad\text{for all }y\in B_{\p}\(0,\alpha\).
\end{equation}
Moreover, by \eqref{Eq1}, we obtain
\begin{equation}\label{Lem3Eq12}
-\Delta_{\p}\widetilde{u}_\alpha=\left|u_\alpha\(y_\alpha\)\right|^{1-p^*}\cdot f_\alpha\(\tau_\alpha\(y\),\left|u_\alpha\(y_\alpha\)\right|\cdot\widetilde{u}_\alpha\)\quad\text{in }\R^n,
\end{equation}
and \eqref{Eq2} gives
\begin{equation}\label{Lem3Eq13}
\left|u_\alpha\(y_\alpha\)\right|^{1-p^*}\cdot\left|f_\alpha\(\tau_\alpha\(y\),\left|u_\alpha\(y_\alpha\)\right|\cdot\widetilde{u}_\alpha\)\right|\le\Lambda\left|\widetilde{u}_\alpha\right|^{p^*-1}.
\end{equation} 
By Lieberman's gradient estimates~\cites{Lie}, it follows from \eqref{Lem3Eq11} and \eqref{Lem3Eq13} that there exists a constant $C>0$ such that for any $R>0$, we have
\begin{equation}\label{Lem3Eq14}
\left\|\nabla\widetilde{u}_\alpha\right\|_{L^\infty(B_{\p}\(0,R\))}\le C
\end{equation}
for large $\alpha$. By Arzela--Ascoli Theorem and a diagonal argument, it follows from \eqref{Lem3Eq11} and \eqref{Lem3Eq14} that $\(\widetilde{u}_\alpha\)_\alpha$ converges up to a subsequence in $C^0_{\loc}\(\R^n\)$ to some Lipschitz continuous function $\widetilde{u}_\infty$ such that $\left|\widetilde{u}_\infty\(0\)\right|=1$. Moreover, by testing \eqref{Lem3Eq12}--\eqref{Lem3Eq13} with $\widetilde{u}_\alpha$, we obtain
\begin{equation}\label{Lem3Eq15}
\sum_{i=1}^n\int_{\R^n}\left|\partial_{x_i}\widetilde{u}_\alpha\right|^{p_i}dx\le\Lambda\int_{\R^n}\left|\widetilde{u}_\alpha\right|^{p^*}dx=\Lambda\int_{\R^n}\left|u_\alpha\right|^{p^*}dx\le\Lambda\(\kappa'\)^{p^*}.
\end{equation}
Since $\left|\partial_{x_i}\left|\widetilde{u}_\alpha\right|\right|=\left|\partial_{x_i}\widetilde{u}_\alpha\right|$ a.e. in $\R^n$, it follows from \eqref{Lem3Eq15} that $\(\left|\widetilde{u}_\alpha\right|\)_\alpha$ converges weakly up to a subsequence to $\left|\widetilde{u}_\infty\right|$ in $D^{1,\p}\(\R^n\)$. Passing to the limit into \eqref{Lem3Eq12}--\eqref{Lem3Eq13}, we then obtain that $\left|\widetilde{u}_\infty\right|$ is a weak solution of the inequality
\begin{equation}\label{Lem3Eq16}
-\Delta_{\p}\left|\widetilde{u}_\infty\right|\le\Lambda\left|\widetilde{u}_\infty\right|^{p^*-1}\quad\text{in }\R^n.
\end{equation}
In particular, since $\left|\widetilde{u}_\infty\(0\)\right|=1$, it follows from Lemma~\ref{Lem2} that $\left\|u_\infty\right\|_{L^{p^*}\(\R^n\)}\ge\kappa_0$, and hence there exists a real number $R>0$ such that 
\begin{equation}\label{Lem3Eq17}
\left\|u_\infty\right\|_{L^{p^*}\(B\(0,R\)\)}>\kappa\,.
\end{equation}
On the other hand, we have
\begin{equation}\label{Lem3Eq18}
\left\|\widetilde{u}_\alpha\right\|_{L^{p^*}(B_{\p}\(0,R\))}=\left\|u_\alpha\right\|_{L^{p^*}\(B_{\p}\(y_\alpha,R\cdot\left|u_\alpha\(y_\alpha\)\right|^{-\delta}\)\)}.
\end{equation}
By \eqref{Lem3Eq8} and since $r_\kappa\(u_\alpha\)<r''$, we obtain
\begin{equation}\label{Lem3Eq19}
B_{\p}\big(y_\alpha,R\cdot\left|u_\alpha\(y_\alpha\)\right|^{-\delta}\big)\cap B_{\p}\(0,r_\kappa\(u_\alpha\)\)=\emptyset
\end{equation}
for large $\alpha$. By definition of $r_\kappa\(u_\alpha\)$, it follows from \eqref{Lem3Eq18} and \eqref{Lem3Eq19} that
\begin{equation}\label{Lem3Eq20}
\left\|\widetilde{u}_\alpha\right\|_{L^{p^*}(B_{\p}\(0,R\))}\le\kappa
\end{equation}
for large $\alpha$, which is in contradiction with \eqref{Lem3Eq17}. This ends the proof of Lemma~\ref{Lem3}.
\endproof

\section{The reverse H\"older-type inequalities}\label{Sec4}

The following result is a key step in the Moser-type iteration scheme that we develop in the next section.

\begin{lemma}\label{Lem4}
Assume that $p_+<p^*$. Let $f:\R^n\times\R\to\R$ be a Caratheodory function such that \eqref{Eq2} holds true, $u$ be a solution of \eqref{Eq1}, and $\kappa$, $r$, and $K_0$ be as in Lemma~\ref{Lem3}. Let $\mathcal{I}_1$ and $\mathcal{I}_2$ be two disjoint subsets of $\left\{1,\dotsc,n\right\}$, $\mathcal{I}_2\ne\emptyset$, and $\q=\(q_i\)_{i\in\mathcal{I}_1\cup\mathcal{I}_2}$ be such that $q_i>1$ for all $i\in\mathcal{I}_1\cup\mathcal{I}_2$. Then there exists a constant $c_0=c_0\(n,\p,\Lambda,K_0\)>1$ such that for any $R_1,R_2>0$, $\lambda<\lambda'\in\(0,1/2\]$, and $\gamma>p_*-1$ such that $\Omega_{\q}\(\mathcal{I}_1,R_1,\mathcal{I}_2,R_2,\lambda'\)\cap B_{\p}\(0,\max\(r,1\)\)=\emptyset$ and $\Omega_{\q}\(\mathcal{I}_1,R_1,\mathcal{I}_2,R_2,\lambda'\)\cap\supp\(u\)$ is bounded, where $\Omega_{\q}\(\mathcal{I}_1,R_1,\mathcal{I}_2,R_2,\lambda\)$ is as in \eqref{Eq6}, we have
\begin{multline}\label{Lem4Eq1}
\left\|u\right\|_{L^\gamma\(\Omega_{\q}\(\mathcal{I}_1,R_1,\mathcal{I}_2,R_2,\lambda\)\)}^\gamma\le c_0\gamma^{p^*}\max_{i\in\mathcal{I}_1\cup\mathcal{I}_2}\Big(\min\(1,\gamma-p_*+1\)^{-p_i}\(\lambda'-\lambda\)^{-p_i}q_i^{p_i}\\
\times R_{\delta_i}^{-\frac{p_i}{q_i}}\left\|u\right\|_{L^{\gamma_i}\(\Omega_{\q}\(\mathcal{I}_1,R_1,\mathcal{I}_2,R_2,\lambda'\)\)}^{\gamma_i}\Big)^\frac{n}{n-p},
\end{multline}
where $\delta_i:=1$ if $i\in\mathcal{I}_1$, $\delta_i:=2$ if $i\in\mathcal{I}_2$, and $\gamma_i:=\frac{n-p}{n}\gamma+p_i-p$ for all $i\in\mathcal{I}_1\cup\mathcal{I}_2$.
\end{lemma}

Preliminary to the proof of Lemma~\ref{Lem4}, we prove the following result.

\begin{lemma}\label{Lem5}
Let $v$ be a nonnegative solution in $D^{1,\p}\(\R^n\)$ of
\begin{equation}\label{Lem5Eq1}
-\Delta_{\p}v\le\Lambda v^{p^*-1}\quad\text{in }\R^n,
\end{equation}
for some real number $\Lambda>0$, where the inequality must be understood in the weak sense as in \eqref{Eq10}. Let $\beta>-1$ and $\eta\in C^1\(\R^n\)$ be such that $0\le\eta\le1$ in $\R^n$, $\eta v$ has compact support, and $\eta^{\(\beta+p_-\)/p_+}\in C^1\(\R^n\)$, where $p_-:=\min\(\left\{p_i\in\p\right\}\)$ and $p_+:=\max\(\left\{p_i\in\p\right\}\)$. Then there exists a constant $C=C\(n,\p\)$ such that
\begin{multline}\label{Lem5Eq2}
\(\int_{\R^n}\(\eta v\)^{\frac{n\(\beta+p\)}{n-p}}dx\)^{\frac{n-p}{n}}\le C\(\left|\beta\right|^p+1\)\bigg(\Lambda\(\beta+1\)^{-1}\int_{\R^n}\eta^{\beta+p_-}v^{\beta+p^*}dx\\
+\sum_{i=1}^n\min\(1,\beta+1\)^{-p_i}\int_{\R^n}\left|\partial_{x_i}\eta\right|^{p_i}\eta^{\beta+p_--p_i}v^{\beta+p_i}dx\bigg).
\end{multline}
\end{lemma}

The finiteness of the integrals in \eqref{Lem5Eq2} is ensured by the fact that $v\in L^\infty\(\R^n\)$, $\eta v$ has compact support, and $\eta^{\(\beta+p_-\)/p_+}\in C^1\(\R^n\)$.

\proof[Proof of Lemma~\ref{Lem5}]
For any $\varepsilon>0$, we define $v_\varepsilon:=v+\varepsilon\overline\eta$, where $\overline\eta$ is a cutoff function on a neighborhood of the support of $\eta v$ such that $\overline\eta^{\beta+1}\in C^1\(\R^n\)$. Since $v\in D^{1,\p}\(\R^n\)\cap L^\infty\(\R^n\)$, $\eta^{\(\beta+p_-\)/p_+}\in C^1\(\R^n\)$, and $\(\beta+p_-\)/p_+\le\(\beta+p_-\)/p_-=1+\beta/p_-$, we get $\(\eta v_\varepsilon\)^{\min\(1,1+\beta/p_-\)}\in D^{1,\p}\(\R^n\)$. By a generalized version of the anisotropic Sobolev inequality (see C\^irstea--V\'etois~\cite{CirVet}*{Lemma~A.1}), we then obtain
\begin{equation}\label{Lem5Eq3}
\left\|\eta v_\varepsilon\right\|_{L^{\frac{n\(\beta+p\)}{n-p}}\(\R^n\)}^{\beta+p}\le C\(\beta+p\)^p\prod_{i=1}^n\left\|\(\eta v_\varepsilon\)^{\frac{\beta}{p_i}}\partial_{x_i}\(\eta v_\varepsilon\)\right\|_{L^{p_i}\(\R^n\)}^{\frac{p}{n}}<\infty
\end{equation}
for some constant $C=C\(n,\p\)$. For any $i=1,\dotsc,n$, we have
\begin{equation}\label{Lem5Eq4}
\int_{\R^n}\(\eta v_\varepsilon\)^\beta\left|\partial_{x_i}\(\eta v_\varepsilon\)\right|^{p_i}dx\le2^{p_i-1}\(\int_{\R^n}\left|\partial_{x_i}\eta\right|^{p_i}\eta^\beta v_\varepsilon^{\beta+p_i}dx+\int_{\R^n}\left|\partial_{x_i} v_\varepsilon\right|^{p_i}\eta^{\beta+p_i}v_\varepsilon^\beta dx\).
\end{equation}
Since $v\in D^{1,\p}\(\R^n\)\cap L^\infty\(\R^n\)$, $\overline\eta^{\beta+1}\in C^1\(\R^n\)$, $\eta^{\(\beta+p_-\)/p_+}\in C^1\(\R^n\)$, and $\(\beta+p_-\)/p_+\le\beta+p_-$, we get $\eta^{\beta+p_i}v_\varepsilon^{\beta+1}\in D^{1,\p}\(\R^n\)$. For any $i=1,\dotsc,n$, since $v_\varepsilon\equiv v+\varepsilon$ on the support of $\eta v$, testing \eqref{Lem5Eq1} with $\eta^{\beta+p_i}v_\varepsilon^{\beta+1}$ gives
\begin{multline}\label{Lem5Eq5}
\(\beta+1\)\sum_{j=1}^n\int_{\R^n}\left|\partial_{x_j}v\right|^{p_j}\eta^{\beta+p_i}v_\varepsilon^\beta dx\le\Lambda\int_{\R^n}\eta^{\beta+p_i}v^{p^*-1}v_\varepsilon^{\beta+1}dx\\
-\(\beta+p_i\)\sum_{j=1}^n\int_{\R^n}\left|\partial_{x_j}v\right|^{p_j-2}\(\partial_{x_j}v\)\(\partial_{x_j}\eta\)\eta^{\beta+p_i-1}v_\varepsilon^{\beta+1}dx\,.
\end{multline}
For any $i,j=1,\dotsc,n$, Youngs inequality yields
\begin{multline}\label{Lem5Eq6}
-\(\beta+p_i\)\left|\partial_{x_j}v\right|^{p_j-2}\(\partial_{x_j}v\)\(\partial_{x_j}\eta\)\eta^{\beta+p_i-1}v_\varepsilon^{\beta+1}\\
\le\frac{p_j-1}{p_j}\cdot\(\beta+1\)\left|\partial_{x_j}v\right|^{p_j}\eta^{\beta+p_i}v_\varepsilon^{\beta}+\frac{1}{p_j}\cdot\frac{\(\beta+p_i\)^{p_j}}{\(\beta+1\)^{p_j-1}}\left|\partial_{x_j}\eta\right|^{p_j}\eta^{\beta+p_i-p_j}v_\varepsilon^{\beta+p_j}.
\end{multline}
It follows from \eqref{Lem5Eq5} and \eqref{Lem5Eq6} that
\begin{multline}\label{Lem5Eq7}
\sum_{j=1}^n\frac{1}{p_j}\int_{\R^n}\left|\partial_{x_j}v\right|^{p_j}\eta^{\beta+p_i}v_\varepsilon^\beta dx\le\Lambda\(\beta+1\)^{-1}\int_{\R^n}\eta^{\beta+p_i}v^{p^*-1}v_\varepsilon^{\beta+1}dx\\
+\sum_{j=1}^n\frac{1}{p_j}\cdot\(\frac{\beta+p_i}{\beta+1}\)^{p_j}\int_{\R^n}\left|\partial_{x_j}\eta\right|^{p_j}\eta^{\beta+p_i-p_j}v_\varepsilon^{\beta+p_j}dx\,.
\end{multline}
In particular, by \eqref{Lem5Eq4} and \eqref{Lem5Eq7}, we obtain
\begin{multline}\label{Lem5Eq8}
\int_{\R^n}\(\eta v_\varepsilon\)^\beta\left|\partial_{x_i}\(\eta v_\varepsilon\)\right|^{p_i}dx\le C\bigg(
\Lambda\(\beta+1\)^{-1}\int_{\R^n}\eta^{\beta+p_i}v^{p^*-1}v_\varepsilon^{\beta+1}dx\\
+\sum_{j=1}^n\min\(1,\beta+1\)^{-p_j}\int_{\R^n}\left|\partial_{x_j}\eta\right|^{p_j}\eta^{\beta+p_i-p_j}v_\varepsilon^{\beta+p_j}dx+\varepsilon^{\beta+p_i}\int_{\R^n}\left|\partial_{x_i}\overline\eta\right|^{p_i}\eta^{\beta+p_i}\overline\eta^{\beta}dx\bigg)
\end{multline}
for some constant $C=C\(n,\p\)$. Finally, since $\eta^{p_i}\le\eta^{p_-}$, we get \eqref{Lem5Eq2} by plugging \eqref{Lem5Eq8} into \eqref{Lem5Eq3} and passing to the limit as $\varepsilon\to0$. This ends the proof of Lemma~\ref{Lem5}.
\endproof

Now, we can prove Lemma~\ref{Lem4} by using Lemma~\ref{Lem5}.

\proof[Proof of Lemma~\ref{Lem4}]
We denote $\beta:=\frac{n-p}{n}\gamma-p$. In particular, $\gamma>p_*-1$ is equivalent to $\beta>-1$. In connexion with the sets $\Omega_{\q}\(\mathcal{I}_1,R_1,\mathcal{I}_2,R_2,\lambda\)$, we define test functions of the form
\begin{equation}\label{Lem4Eq2}
\eta\(x\):=\bigg[\overline\eta_{\lambda,\lambda'}\bigg(R_1^{-1}\sum_{i\in\mathcal{I}_1}\left|x_i\right|^{q_i}\bigg)\widetilde\eta_{\lambda,\lambda'}\bigg(R_2^{-1}\sum_{i\in\mathcal{I}_2}\left|x_i\right|^{q_i}\bigg)\bigg]^{\max\(1,\frac{p_+}{\beta+p_-}\)}
\end{equation}
for all $x\in\R^n$, where $\overline\eta_{\lambda,\lambda'},\,\widetilde\eta_{\lambda,\lambda'}\in C^1\(0,\infty\)$ satisfy $0\le\overline\eta_{\lambda,\lambda'},\,\widetilde\eta_{\lambda,\lambda'}\le1$ in $\(0,\infty\)$, $\overline\eta_{\lambda,\lambda'}=1$ in $\[0,1+\lambda\]$, $\overline\eta_{\lambda,\lambda'}=0$ in $\[1+\lambda',\infty\)$, $\left|\overline\eta'_{\lambda,\lambda'}\right|\le2$ in $\[1+\lambda,1+\lambda'\]$, $\widetilde\eta_{\lambda,\lambda'}=1$ in $\[1-\lambda,1+\lambda\]$, $\widetilde\eta_{\lambda,\lambda'}=0$ in $\[0,1-\lambda'\]\cup\[1+\lambda',\infty\)$, and $\left|\widetilde\eta'_{\lambda,\lambda'}\right|\le2/\(\lambda'-\lambda\)$ in $\[1-\lambda',1-\lambda\]\cup\[1+\lambda,1+\lambda'\]$. With these properties of $\overline\eta_{\lambda,\lambda'}$ and $\widetilde\eta_{\lambda,\lambda'}$, we obtain
$$0\le\eta\le1\text{ in }\R^n,\,\eta=1\text{ in }\Omega_{\q}\(\mathcal{I}_1,R_1,\mathcal{I}_2,R_2,\lambda\),\text{ and }\eta=0\text{ in }\R^n\backslash \Omega_{\q}\(\mathcal{I}_1,R_1,\mathcal{I}_2,R_2,\lambda'\).$$
Since $\Omega_{\q}\(\mathcal{I}_1,R_1,\mathcal{I}_2,R_2,\lambda'\)\cap\supp\(u\)$ is bounded by assumption, we get that $\eta u$ has compact support. Moreover, since $q_i>1$ for all $i\in\mathcal{I}_1\cup\mathcal{I}_2$, we get $\eta^{\(\beta+p_-\)/p_+}\in C^1\(\R^n\)$ and
\begin{equation}\label{Lem4Eq3}
\left|\partial_{x_i}\eta\(x\)\right|^{p_i}\eta\(x\)^{\beta+p_--p_i}\le\(\frac{4q_i}{\lambda'-\lambda}\max\(1,\frac{p_+}{\beta+p_-}\)\)^{p_i}R_{\delta_i}^{-\frac{p_i}{q_i}}
\end{equation}
for all $x\in\supp\(\eta\)$, where $\delta_i:=1$ if $i\in\mathcal{I}_1$, $\delta_i:=2$ if $i\in\mathcal{I}_2$. By applying Lemma~\ref{Lem5} with $v=\left|u\right|$ and $\eta$ as in \eqref{Lem4Eq2}, and using \eqref{Lem4Eq3}, we obtain
\begin{multline}\label{Lem4Eq4}
\(\int_{\Omega_{\q}\(\mathcal{I}_1,R_1,\mathcal{I}_2,R_2,\lambda\)}u^\gamma dx\)^{\frac{n-p}{n}}\le C\(\left|\beta\right|^p+1\)\bigg(\Lambda\(\beta+1\)^{-1}\int_{\Omega_{\q}\(\mathcal{I}_1,R_1,\mathcal{I}_2,R_2,\lambda'\)}u^{\beta+p^*}dx\\
+\sum_{i\in\mathcal{I}_1\cup\mathcal{I}_2}\min\(1,\beta+1\)^{-p_i}\(\lambda'-\lambda\)^{-p_i}q_i^{p_i}R_{\delta_i}^{-\frac{p_i}{q_i}}\int_{\Omega_{\q}\(\mathcal{I}_1,R_1,\mathcal{I}_2,R_2,\lambda'\)}u^{\beta+p_i}dx\bigg)
\end{multline}
for some constant $C=C\(n,\p\)$. 

Now, we estimate the first integral in the right-hand side of \eqref{Lem4Eq4}. We claim that there exists a constant $C'=C\(n,\p,K_0\)$ such that 
\begin{equation}\label{Lem4Eq5}
u\(x\)^{p^*-p_{i_0}}\le C'R_2^{-\frac{p_{i_0}}{q_{i_0}}}\quad\text{for all }x\in\Omega_{\q}\(\mathcal{I}_1,R_1,\mathcal{I}_2,R_2,\lambda'\),
\end{equation}
where $K_0$ is the constant given by Lemma~\ref{Lem3} and $i_0\in\mathcal{I}_2$ is such that
\begin{equation}\label{Lem4Eq6}
\frac{q_{i_0}\(p^*-p_{i_0}\)}{p_{i_0}}=\max_{i\in\mathcal{I}_2}\(\frac{q_i\(p^*-p_i\)}{p_i}\)\,.
\end{equation}
We prove this claim. For any $x\in\Omega_{\q}\(\mathcal{I}_1,R_1,\mathcal{I}_2,R_2,\lambda'\)$, since $\lambda'\le1/2$, we obtain
\begin{equation}\label{Lem4Eq7}
\frac{R_2}{2}\le\sum_{i\in\mathcal{I}_2}\left|x_i\right|^{q_i}\le n\cdot|x_{i\(x\)}|^{q_{i\(x\)}}\le n\cdot d_{\p}\(x,0\)^{\frac{q_{i\(x\)}(p^*-p_{i\(x\)})}{\delta p_{i\(x\)}}},
\end{equation}
where $i\(x\)\in\mathcal{I}_2$ is such that $|x_{i\(x\)}|^{q_{i\(x\)}}=\max\(\left\{\left|x_i\right|^{q_i}:\,i\in\mathcal{I}_2\right\}\)$, and the distance function $d_{\p}$ and the real number $\delta$ are as in \eqref{Lem3Eq2}. Since $\Omega_{\q}\(\mathcal{I}_1,R_1,\mathcal{I}_2,R_2,\lambda'\)\cap B_{\p}\(0,\max\(r,1\)\)=\emptyset$ by assumption, \eqref{Lem4Eq5} follows from \eqref{Lem4Eq6}, \eqref{Lem4Eq7}, and Lemma~\ref{Lem3}. In particular, \eqref{Lem4Eq5} implies
\begin{equation}\label{Lem4Eq9}
\int_{\Omega_{\q}\(\mathcal{I}_1,R_1,\mathcal{I}_2,R_2,\lambda'\)}u^{\beta+p^*}dx\le C'R_2^{-\frac{p_{i_0}}{q_{i_0}}}\int_{\Omega_{\q}\(\mathcal{I}_1,R_1,\mathcal{I}_2,R_2,\lambda'\)}u^{\beta+p_{i_0}}dx\,.
\end{equation}
Finally, \eqref{Lem4Eq1} follows from \eqref{Lem4Eq4}, \eqref{Lem4Eq9}, and the fact that $\beta+1=\frac{n-p}{n}\(\gamma-p_*+1\)$ and $\beta+p_i=\gamma_i$. This ends the proof of Lemma~\ref{Lem4}.
\endproof

\section{The iteration scheme}\label{Sec5}

In this section, we describe the iteration scheme which leads to the proofs of our main results.

\medskip
Let $\mathcal{I}_1$ and $\mathcal{I}_2$ be two disjoint subsets of $\left\{1,\dotsc,n\right\}$, $\mathcal{I}_2\ne\emptyset$, and $\q=\(q_i\)_{i\in\mathcal{I}_1\cup\mathcal{I}_2}$ be such that $q_i>1$ for all $i\in\mathcal{I}_1\cup\mathcal{I}_2$. The idea is to apply Lemma~\ref{Lem4} by induction. For any $\gamma>p_*-1$, Lemma~\ref{Lem4} provides an estimate of the $L^\gamma$--norm of $u$ with respect to the set of $L^{\gamma_{i_1}}$--norms of $u$, where $\gamma_{i_1}:=\frac{n-p}{n}\gamma+p_{i_1}-p$ for all $i_1\in\mathcal{I}_1\cup\mathcal{I}_2$. If $\gamma_{i_1}>p_*-1$, then another application of Lemma~\ref{Lem4} gives estimates of the $L^{\gamma_{i_1}}$--norms of $u$ with respect to the set of $L^{\gamma_{i_1i_2}}$--norms of $u$, where $\gamma_{i_1i_2}:=\frac{n-p}{n}\gamma_{i_1}+p_{i_2}-p$, etc... By induction, we define
\begin{equation}\label{Eq11}
\gamma_{i_1,\dotsc,i_{j+1}}:=\frac{n-p}{n}\gamma_{i_1,\dotsc,i_j}+p_{i_{j+1}}-p
\end{equation}
for all $j\in\N$ and $i_1,\dotsc,i_{j+1}\in\mathcal{I}_1\cup\mathcal{I}_2$, with the convention that $\gamma_{i_1,\dotsc,i_j}:=\gamma$ if $j=0$. In particular, we obtain the formula
\begin{equation}\label{Eq12}
\gamma_{i_1,\dotsc,i_k}=\(\frac{n-p}{n}\)^k\gamma+\sum_{j=1}^k\(\frac{n-p}{n}\)^{k-j}\(p_{i_j}-p\)
\end{equation}
for all $k\in\N$. The stopping condition in our induction argument is $\gamma_{i_1,\dotsc,i_k}<\frac{n}{p}\(p_\varepsilon-p\)$, where
\begin{equation}\label{Eq13}
p_\varepsilon:=\(1+\varepsilon\)p_0\,,\quad p_0:=\max\(p_*,\left\{p_i\,:\,i\in\mathcal{I}_1\cup\mathcal{I}_2\right\}\),
\end{equation}
and $\varepsilon$ is a fixed real number in $\(0,1\)$. Note that $\frac{n}{p}\(p_\varepsilon-p\)>\frac{n}{p}\(p_*-p\)=p_*-1$ so that we can apply Lemma~\ref{Lem4} as long as our stopping condition is not satisfied. For any $k\ge1$, we let $\Phi_{k,\gamma,\varepsilon}$ be the set of all sequences of indices for which our induction argument stops after exactly $k$ iterations, namely
\begin{multline}\label{Eq14}
\Phi_{k,\gamma,\varepsilon}:=\Big\{\(i_1,\dotsc,i_k\)\in\(\mathcal{I}_1\cup\mathcal{I}_2\)^k\,:\quad\gamma_{i_1,\dotsc,i_j}\ge\frac{n}{p}\(p_\varepsilon-p\)\text{ for all }j=0,\dotsc,k-1\\
\text{and}\quad\gamma_{i_1,\dotsc,i_k}<\frac{n}{p}\(p_\varepsilon-p\)\Big\}.
\end{multline}
The following result provides a control on the number of iterations in our induction argument.

\begin{lemma}\label{Lem6}
Let $\mathcal{I}_1$ and $\mathcal{I}_2$ be two disjoint subsets of $\left\{1,\dotsc,n\right\}$, $\mathcal{I}_2\ne\emptyset$, and $\q=\(q_i\)_{i\in\mathcal{I}_1\cup\mathcal{I}_2}$ be such that $q_i>1$ for all $i\in\mathcal{I}_1\cup\mathcal{I}_2$. Then for any $\varepsilon>0$, $\gamma\ge\frac{n}{p}\(p_\varepsilon-p\)$, $k\in\N$, and $\(i_1,\dotsc,i_k\)\in\(\mathcal{I}_1\cup\mathcal{I}_2\)^k$, we have
\begin{equation}\label{Lem6Eq1}
\gamma_{i_1,\dotsc,i_k}>\frac{n}{p}\(p_\varepsilon-p\)\quad\text{if }k<k_{\gamma,\varepsilon}^-\quad\text{and}\quad\gamma_{i_1,\dotsc,i_k}<\frac{n}{p}\(p_\varepsilon-p\)\quad\text{if }k\ge k_{\gamma,\varepsilon}^+\,,
\end{equation}
where $\gamma_{i_1,\dotsc,i_k}$ is as in \eqref{Eq11}, $p_\varepsilon$ is as in \eqref{Eq13}, and $k_{\gamma,\varepsilon}^-$ and $k_{\gamma,\varepsilon}^+$ are the smallest and largest natural numbers, respectively, such that 
\begin{equation}\label{Lem6Eq2}
\frac{n}{p}\(\frac{n}{n-p}\)^{k_{\gamma,\varepsilon}^+-1}\varepsilon p_0<\gamma<\frac{n}{p}\(\frac{n}{n-p}\)^{k_{\gamma,\varepsilon}^-}\(p_\varepsilon-p_-\),
\end{equation}
where $p_-:=\min\(\left\{p_i\in\p\right\}\)$. In particular, we have $\Phi_{k,\gamma,\varepsilon}=\emptyset$ for all $k<k_{\gamma,\varepsilon}^-$ and $k>k_{\gamma,\varepsilon}^+$, where $\Phi_{k,\gamma,\varepsilon}$ is as in \eqref{Eq14}.
\end{lemma}

\proof[Proof of Lemma~\ref{Lem6}]
Since $p_-\le p_{i_j}\le p_0$ for all $j=1,\dotsc,k$, it follows from \eqref{Eq12} that
\begin{equation}\label{Lem6Eq3}
-\sum_{j=1}^k\(\frac{n-p}{n}\)^{k-j}\(p-p_-\)\le\gamma_{i_1,\dotsc,i_k}-\(\frac{n-p}{n}\)^k\gamma\le\sum_{j=1}^k\(\frac{n-p}{n}\)^{k-j}\(p_0-p\).
\end{equation}
Moreover, by a simple calculation, we obtain
\begin{equation}\label{Lem6Eq4}
\sum_{j=1}^k\(\frac{n-p}{n}\)^{k-j}=\frac{n}{p}\(1-\(\frac{n-p}{n}\)^k\)<\frac{n}{p}\,.
\end{equation}
It follows from \eqref{Lem6Eq3} and \eqref{Lem6Eq4} that
\begin{equation}\label{Lem6Eq5}
-\frac{n}{p}\(p-p_-\)<\gamma_{i_1,\dotsc,i_k}-\(\frac{n-p}{n}\)^k\gamma<\frac{n}{p}\(p_0-p\).
\end{equation}
Finally, \eqref{Lem6Eq1} follows from \eqref{Lem6Eq5} together with the definitions of $k_{\gamma,\varepsilon}^-$ and $k_{\gamma,\varepsilon}^+$.
\endproof 

Now, we can prove the main result of this section.

\begin{lemma}\label{Lem7}
Assume that $p_+<p^*$. Let $f:\R^n\times\R\to\R$ be a Caratheodory function such that \eqref{Eq2} holds true, $u$ be a solution of \eqref{Eq1}, and $\kappa$, $r$, and $K_0$ be as in Lemma~\ref{Lem3}. Let $\mathcal{I}_1$ and $\mathcal{I}_2$ be two disjoint subsets of $\left\{1,\dotsc,n\right\}$, $\mathcal{I}_2\ne\emptyset$, and $\q=\(q_i\)_{i\in\mathcal{I}_1\cup\mathcal{I}_2}$ be such that $q_i>1$ for all $i\in\mathcal{I}_1\cup\mathcal{I}_2$. Then there exists a constant $c_1=c_1\(n,\p,\Lambda,K_0\)>1$ such that for any $\varepsilon\in\(0,1\)$, $\gamma>\frac{n}{p}\(p_\varepsilon-p\)$, and $R_1,R_2>0$ such that $\Omega_{\q}\(\mathcal{I}_1,R_1,\mathcal{I}_2,R_2,1/2\)\cap B_{\p}\(0,\max\(r,1\)\)=\emptyset$ and $\Omega_{\q}\(\mathcal{I}_1,R_1,\mathcal{I}_2,R_2,1/2\)\cap\supp\(u\)$ is bounded, we have
\begin{multline}\label{Lem7Eq1}
\left\|u\right\|_{L^{\gamma}\(\Omega_{\q}\(\mathcal{I}_1,R_1,\mathcal{I}_2,R_2,\lambda_{0,\gamma,\varepsilon}\)\)}\le c_1^{\frac{1}{\varepsilon}}\max\(q_i\)^{\frac{1}{\varepsilon}\cdot\frac{n}{n-p}}\(p_\varepsilon-p_*\)^{-\frac{1}{\varepsilon}\cdot\frac{n}{n-p}}\\
\times\max_{\(i_1,\dotsc,i_k\)\in\Phi_{\gamma,\varepsilon}}\(\(\prod_{j=1}^kR_{\delta_{i_j}}^{-\frac{1}{\gamma}\(\frac{n}{n-p}\)^j\frac{p_{i_j}}{q_{i_j}}}\)\left\|u\right\|_{L^{\gamma_{i_1,\dotsc,i_k}}\(\Omega_{\q}\(\mathcal{I}_1,R_1,\mathcal{I}_2,R_2,\lambda_{k,\gamma,\varepsilon}\)\)}^{\frac{\gamma_{i_1,\dotsc,i_k}}{\gamma}\(\frac{n}{n-p}\)^k}\),
\end{multline}
where $\delta_{i_j}:=1$ if $i_j\in\mathcal{I}_1$, $\delta_{i_j}:=2$ if $i_j\in\mathcal{I}_2$, $\gamma_{i_1,\dotsc,i_k}$ is as in \eqref{Eq11}, $p_\varepsilon$ is as in \eqref{Eq13}, $\Omega_{\q}\(\mathcal{I}_1,R_1,\mathcal{I}_2,R_2,\lambda_k\)$ is as in \eqref{Eq6}, and
\begin{equation}\label{Lem7Eq2}
\lambda_{k,\gamma,\varepsilon}:=\frac{1}{4}\big(1+2^{k-k_{\gamma,\varepsilon}^+-1}\big)\quad\text{and}\quad\Phi_{\gamma,\varepsilon}:=\bigcup_{k=k_{\gamma,\varepsilon}^-}^{k_{\gamma,\varepsilon}^+}\Phi_{k,\gamma,\varepsilon}
\end{equation}
with $k_{\gamma,\varepsilon}^-$ and $k_{\gamma,\varepsilon}^+$ as in Lemma~\ref{Lem6}, and $\Phi_{k,\gamma,\varepsilon}$ as in \eqref{Eq14}.
\end{lemma}

\proof[Proof of Lemma~\ref{Lem7}]
Applying Lemma~\ref{Lem4} by induction with the stopping condition $\gamma_{i_1,\dotsc,i_k}<\frac{n}{p}\(p_\varepsilon-p\)$ gives
\begin{multline}\label{Lem7Eq3}
\left\|u\right\|_{L^{\gamma}\(\Omega_{\q}\(\mathcal{I}_1,R_1,\mathcal{I}_2,R_2,\lambda_{0,\gamma,\varepsilon}\)\)}\le\max_{\(i_1,\dotsc,i_k\)\in\Phi_{\gamma,\varepsilon}}\Big(\mathcal{A}_{k,\gamma}\times\mathcal{B}_{i_1,\dotsc,i_{k-1},\gamma}\times\mathcal{C}_{i_1,\dotsc,i_k,\gamma}\times\mathcal{D}_{i_1,\dotsc,i_k,\gamma,\varepsilon}\\
\times\(\prod_{j=1}^kR_{\delta_{i_j}}^{-\frac{1}{\gamma}\(\frac{n}{n-p}\)^j\frac{p_{i_j}}{q_{i_j}}}\)\left\|u\right\|_{L^{\gamma_{i_1,\dotsc,i_k}}\(\Omega_{\q}\(\mathcal{I}_1,R_1,\mathcal{I}_2,R_2,\lambda_{k,\gamma,\varepsilon}\)\)}^{\frac{\gamma_{i_1,\dotsc,i_k}}{\gamma}\(\frac{n}{n-p}\)^k}\Big),
\end{multline}
where $\lambda_{k,\gamma,\varepsilon}$ and $\Phi_{\gamma,\varepsilon}$ are as in \eqref{Lem7Eq2}, and
\begin{align*}
&\mathcal{A}_{k,\gamma}:=\big(c_0\cdot\max\big(q_i^{\frac{np_i}{n-p}}\big)\big)^{\frac{1}{\gamma}\overset{k-1}{\underset{j=0}{\sum}}\(\frac{n}{n-p}\)^j},\quad\mathcal{B}_{i_1,\dotsc,i_{k-1},\gamma}:=\prod_{j=1}^k\gamma_{i_1,\dotsc,i_{j-1}}^{\frac{p}{\gamma}\(\frac{n}{n-p}\)^j},\allowdisplaybreaks\\
&\mathcal{C}_{i_1,\dotsc,i_{k-1},\gamma}:=\prod_{j=1}^k\min\(1,\gamma_{i_1,\dotsc,i_{j-1}}-p_*+1\)^{-\frac{p_{i_j}}{\gamma}\(\frac{n}{n-p}\)^j},\allowdisplaybreaks\\
&\mathcal{D}_{i_1,\dotsc,i_k,\gamma,\varepsilon}:=\prod_{j=1}^k\(\lambda_{j,\gamma,\varepsilon}-\lambda_{j-1,\gamma,\varepsilon}\)^{-\frac{p_{i_j}}{\gamma}\(\frac{n}{n-p}\)^j}.
\end{align*}
Now, we fix $\(i_1,\dotsc,i_k\)\in\Phi_{\gamma,\varepsilon}$ and estimate each of the terms in the right-hand side of \eqref{Lem7Eq3}. 

\proof[Estimate of $\mathcal{A}_{k,\gamma}$]
By using the fact that $k\le k_{\gamma,\varepsilon}^+$ and applying \eqref{Lem6Eq2}, we obtain
\begin{equation}\label{Lem7Eq4}
\sum_{j=0}^{k-1}\(\frac{n}{n-p}\)^j=\frac{n-p}{p}\[\(\frac{n}{n-p}\)^k-1\]<\frac{n}{p}\(\frac{n}{n-p}\)^{k_{\gamma,\varepsilon}^+-1}<\frac{\gamma}{\varepsilon p_0}\,.
\end{equation}
Since $c_0>1$, $q_i>1$, and $p_i\le p_0$ for all $i\in\mathcal{I}_1\cup\mathcal{I}_2$, it follows from \eqref{Lem7Eq4} that
\begin{equation}\label{Lem7Eq5}
\mathcal{A}_{k,\gamma}<c_0^{\frac{1}{\varepsilon p_0}}\max\(q_i\)^{\frac{1}{\varepsilon}\cdot\frac{n}{n-p}}.
\end{equation}

\proof[Estimate of $\mathcal{B}_{i_1,\dotsc,i_{k-1},\gamma}$]
For any $j=1,\dotsc,k$, since $p_i\le p_0$ for all $i\in\mathcal{I}_1\cup\mathcal{I}_2$, by \eqref{Eq12}, \eqref{Lem6Eq2}, and \eqref{Lem6Eq4},  we obtain
\begin{equation}\label{Lem7Eq6}
\gamma_{i_1,\dotsc,i_{j-1}}\le\(\frac{n-p}{n}\)^{j-1}\gamma+\frac{n}{p}\(p_0-p\)\le C\max\bigg(1,\(\frac{n}{n-p}\)^{k_{\gamma,\varepsilon}^--j}\bigg)
\end{equation}
for some constant $C=C\(n,\p\)>1$. It follows from \eqref{Lem7Eq6} that
\begin{equation}\label{Lem7Eq7}
\mathcal{B}_{i_1,\dotsc,i_{k-1},\gamma}\le C^{\frac{p}{\gamma}\overset{k}{\underset{j=1}{\sum}}\(\frac{n}{n-p}\)^j}\(\frac{n}{n-p}\)^{\frac{p}{\gamma}\overset{k_{\gamma,\varepsilon}^-}{\underset{j=1}{\sum}}\(\frac{n}{n-p}\)^j\(k_{\gamma,\varepsilon}^--j\)}.
\end{equation}
A simple calculation gives
$$\overset{k_{\gamma,\varepsilon}^-}{\underset{j=1}{\sum}}\(\frac{n}{n-p}\)^j\(k_{\gamma,\varepsilon}^--j\)=\frac{n^2}{p^2}\[\(\frac{n}{n-p}\)^{k_{\gamma,\varepsilon}^--1}-\frac{p}{n}\(k_{\gamma,\varepsilon}^--1\)-1\]<\frac{n^2}{p^2}\(\frac{n}{n-p}\)^{k_{\gamma,\varepsilon}^--1},$$
and hence by definition of $k_{\gamma,\varepsilon}^-$, we obtain
\begin{equation}\label{Lem7Eq8}
\overset{k_{\gamma,\varepsilon}^-}{\underset{j=1}{\sum}}\(\frac{n}{n-p}\)^j\(k_{\gamma,\varepsilon}^--j\)\le\frac{n}{p}\cdot\frac{\gamma}{p_\varepsilon-p_-}<\frac{n}{p}\cdot\frac{\gamma}{p_*-p_-}\,.
\end{equation}
It follows from \eqref{Lem7Eq4}, \eqref{Lem7Eq7}, and \eqref{Lem7Eq8} that
\begin{equation}\label{Lem7Eq9}
\mathcal{B}_{i_1,\dotsc,i_{k-1},\gamma}\le C^{\frac{p^*}{\varepsilon p_0}}\(\frac{n}{n-p}\)^{\frac{n}{p_*-p_-}}.
\end{equation}

\proof[Estimate of $\mathcal{C}_{i_1,\dotsc,i_{k-1},\gamma}$]
Since $p_*-1=\frac{n}{p}\(p_*-p\)$ and $\gamma_{i_1,\dotsc,i_{j-1}}>\frac{n}{p}\(p_\varepsilon-p\)$ for all $j=1,\dotsc,k$, we obtain
\begin{equation}\label{Lem7Eq10}
\mathcal{C}_{i_1,\dotsc,i_{k-1},\gamma}\le\min\Big(1,\frac{n}{p}\(p_\varepsilon-p_*\)\Big)^{-\frac{1}{\gamma}\overset{k}{\underset{j=1}{\sum}}p_{i_j}\(\frac{n}{n-p}\)^j}.
\end{equation}
Since $p_{i_j}\le p_0$ for all $j=1,\dotsc,k$, it follows from \eqref{Lem7Eq4} and \eqref{Lem7Eq10} that
\begin{equation}\label{Lem7Eq11}
\mathcal{C}_{i_1,\dotsc,i_{k-1},\gamma}\le\min\Big(1,\frac{n}{p}\(p_\varepsilon-p_*\)\Big)^{-\frac{1}{\varepsilon}\cdot\frac{n}{n-p}}.
\end{equation}

\proof[Estimate of $\mathcal{D}_{i_1,\dotsc,i_k,\gamma,\varepsilon}$]
By \eqref{Lem7Eq2} and since $k\le k_{\gamma,\varepsilon}^+$ and $p_{i_j}\le p_0$ for all $j=1,\dotsc,k$, we obtain
\begin{equation}\label{Lem7Eq12}
\mathcal{D}_{i_1,\dotsc,i_k,\gamma,\varepsilon}\le2^{\frac{1}{\gamma}\overset{k}{\underset{j=1}{\sum}}p_{i_j}\(\frac{n}{n-p}\)^j\(k_{\gamma,\varepsilon}^+-j+4\)}\le2^{\frac{p_0}{\gamma}\overset{k_{\gamma,\varepsilon}^+}{\underset{j=1}{\sum}}\(\frac{n}{n-p}\)^j\(k_{\gamma,\varepsilon}^+-j+4\)}.
\end{equation}
We find
\begin{align*}
\overset{k_{\gamma,\varepsilon}^+}{\underset{j=1}{\sum}}\(\frac{n}{n-p}\)^j\(k_{\gamma,\varepsilon}^+-j+4\)&=\frac{n}{p}\[\(\frac{n}{n-p}\)^{k_{\gamma,\varepsilon}^+}\(3+\frac{n}{p}\)-k_{\gamma,\varepsilon}^+-3-\frac{n}{p}\]\\
&<\frac{n}{p}\(\frac{n}{n-p}\)^{k_{\gamma,\varepsilon}^+}\(3+\frac{n}{p}\),
\end{align*}
and hence by \eqref{Lem6Eq2}, we obtain
\begin{equation}\label{Lem7Eq13}
\overset{k_{\gamma,\varepsilon}^+}{\underset{j=1}{\sum}}\(\frac{n}{n-p}\)^j\(k_{\gamma,\varepsilon}^+-j+4\)<\frac{\gamma}{\varepsilon p_0}\cdot\frac{n}{n-p}\(3+\frac{n}{p}\).
\end{equation}
It follows from \eqref{Lem7Eq12} and \eqref{Lem7Eq13} that
\begin{equation}\label{Lem7Eq14}
\mathcal{D}_{i_1,\dotsc,i_k,\gamma,\varepsilon}\le2^{\frac{1}{\varepsilon}\cdot\frac{n}{n-p}\(3+\frac{n}{p}\)}.
\end{equation}

\proof[End of proof of Lemma~\ref{Lem7}]
The estimate \eqref{Lem7Eq1} follows from \eqref{Lem7Eq3}, \eqref{Lem7Eq5}, \eqref{Lem7Eq9}, \eqref{Lem7Eq11}, and  \eqref{Lem7Eq14}.
\endproof

\section{The vanishing result}\label{Sec6}

In this section, we prove a vanishing result which will give Point (i) in Theorem~\ref{Th3}. We define 
\begin{equation}\label{Eq15}
\overline{p}_0:=\max\(p_*,\left\{p_i\in\p\,:\,i\in\Theta\right\}\),
\end{equation}
where $p_*$ is as in \eqref{Eq4} and $\Theta$ is the set of all indices $i$ such that
\begin{equation}\label{Eq16}
\big(p_i-p_--\frac{n}{p}\(p_i-p_*\)\big)\sum_{j=1}^n\max\(\frac{p_i-p_j}{p_j},0\)\ge\(p_*-1\)\(p_i-p_-\)
\end{equation}
with $p_-:=\min\(\left\{p_i\in\p\right\}\)$. We define $\mathcal{I}_0$ as the set of all indices $i$ such that $p_i>\overline{p}_0$.

\medskip
When $p_+>p_*$, one easily sees that the condition \eqref{Eq16} does not hold true for $p_i=p_+$, and hence we have $\overline{p}_0<p_+$ and $\mathcal{I}_0\ne\emptyset$.

\medskip
We prove the following result.

\begin{theorem}\label{Th6}
Assume that $p_*<p_+<p^*$. Let $f:\R^n\times\R\to\R$ be a Caratheodory function such that \eqref{Eq2} holds true and $u$ be a solution of \eqref{Eq1}. Then there exists a constant $R_0=R_0\(n,\p,\Lambda,u\)$ such that $u\(x\)=0$ for all $x\in\R^n$ such that $\sum_{i\in \mathcal{I}_0}\left|x_i\right|\ge R_0$.
\end{theorem}

The proof of Theorem~\ref{Th6} is based on the following result, which we obtain by applying the iteration scheme in Section \ref{Sec5}.

\begin{lemma}\label{Lem8}
Assume that $p_*<p_+<p^*$. Let $f:\R^n\times\R\to\R$ be a Caratheodory function such that \eqref{Eq2} holds true, $u$ be a solution of \eqref{Eq1}, and $\kappa$, $r$, and $K_0$ be as in Lemma~\ref{Lem3}. Let $\overline{p}_0$ be as in \eqref{Eq15} and $p_0\in\p$ be such that
\begin{equation}\label{Lem8Eq1}
p_0>\overline{p}_0\quad\text{and}\quad R_i\(u\)<\infty\text{ for all indices }i\text{ such that }p_i>p_0\,,
\end{equation}
where
\begin{equation}\label{Lem8Eq2}
R_i\(u\):=\sup\(\left\{\left|x_i\right|\,:\,x\in\supp\(u\)\right\}\).
\end{equation}
Let $\mathcal{I}_1$, $\mathcal{I}_2$ be the sets of indices $i$ such that $p_i<p_0$, $p_i=p_0$, respectively. For any $\varepsilon,\lambda\in\(0,1\)$ and $R>1$, we define
\begin{equation}\label{Lem8Eq3}
A_\varepsilon\(R,\lambda\):=\Omega_{\q}\(\mathcal{I}_1,R^{1/\varepsilon},\mathcal{I}_2,R,\lambda\)\quad\text{with}\quad q_i:=\left\{\begin{aligned}&\frac{p_\varepsilon p_i}{p_\varepsilon-p_i}&&\text{if }i\in\mathcal{I}_1\,,\\&p_\varepsilon&&\text{if }i\in\mathcal{I}_2\,,\end{aligned}\right.
\end{equation}
where $p_\varepsilon:=\(1+\varepsilon\)p_0$. If $A_\varepsilon\(R,1/2\)\cap B_{\p}\(0,\max\(r,1\)\)=\emptyset$, then
\begin{equation}\label{Lem8Eq4}
\left\|u\right\|_{L^\infty\(A_\varepsilon\(R,1/4\)\)}\le \big(c_2R^{-\frac{1}{p_\varepsilon}}\big)^{\frac{1}{\varepsilon}}
\end{equation}
for some constant $c_2=c_2(n,\p,\Lambda,K_0,\left\|u\right\|_{L^{p_*-1,\infty}\(\R^n\)},R_0\(u\))$, where
\begin{equation}\label{Lem8Eq5}
R_0\(u\):=\max\(\left\{R_i\(u\):i\in\left\{1,\dotsc,n\right\}\backslash\(\mathcal{I}_1\cup\mathcal{I}_2\)\right\}\).
\end{equation}
\end{lemma}

\proof[Proof of Lemma~\ref{Lem8}]
As is easily seen, we have $1<q_i<s_0$ for some constant $s_0=s_0\(\p\)$. Moreover, by \eqref{Lem8Eq1}, we obtain that $p_\varepsilon-p_*>p_0-p_*>0$ and $A_\varepsilon\(R,1/2\)\cap\supp\(u\)$ is bounded. By Lemma~\ref{Lem7}, we then get that there exists a constant $\widetilde{c}_1=\widetilde{c}_1\(n,\p,\Lambda,K_0\)$ such that for any $\gamma>\frac{n}{p}\(p_\varepsilon-p\)$, we have
\begin{equation}\label{Lem8Eq6}
\left\|u\right\|_{L^{\gamma}\(A_\varepsilon\(R,\lambda_{0,\gamma,\varepsilon}\)\)}\le \widetilde{c}_1^{\,\frac{1}{\varepsilon}}\max_{\(i_1,\dotsc,i_k\)\in\Phi_{\gamma,\varepsilon}}\Big(R^{-\sigma_{i_1,\dotsc,i_k,\gamma,\varepsilon}}\left\|u\right\|_{L^{\gamma_{i_1,\dotsc,i_k}}\(A_\varepsilon\(R,\lambda_{k,\gamma,\varepsilon}\)\)}^{\frac{\gamma_{i_1,\dotsc,i_k}}{\gamma}\(\frac{n}{n-p}\)^k}\Big)
\end{equation}
provided that $A_\varepsilon\(R,1/2\)\cap B_{\p}\(0,\max\(r,1\)\)=\emptyset$, where $\gamma_{i_1,\dotsc,i_k}$ is as in \eqref{Eq11}, $\lambda_{k,\gamma,\varepsilon}$ and $\Phi_{\gamma,\varepsilon}$ are as in \eqref{Lem7Eq2}, and
\begin{equation}\label{Lem8Eq7}
\sigma_{i_1,\dotsc,i_k,\gamma,\varepsilon}:=\frac{1}{\varepsilon\gamma p_\varepsilon}\sum_{j=1}^k\(\frac{n}{n-p}\)^j\(p_\varepsilon-p_{i_j}\)\,.
\end{equation}

We claim that for any $\nu\in\(0,1\)$ and $\(i_1,\dotsc,i_k\)\in\Phi_{\gamma,\varepsilon}$, there exists a constant $c_\nu=c\big(n,\p,K_0,\left\|u\right\|_{L^{p_*-1,\infty}\(\R^n\)},R_0\(u\),\nu\big)$ such that
\begin{equation}\label{Lem8Eq8}
\left\|u\right\|_{L^{\gamma_{i_1,\dotsc,i_k}}\(A_\varepsilon\(R,\lambda_{k,\gamma,\varepsilon}\)\)}^{\frac{\gamma_{i_1,\dotsc,i_k}}{\gamma}\(\frac{n}{n-p}\)^k}\le c_\nu^{\frac{1}{\varepsilon}}R^{\tau_{i_1,\dotsc,i_k,\gamma,\varepsilon,\nu}},
\end{equation}
where
\begin{equation}\label{Lem8Eq9}
\tau_{i_1,\dotsc,i_k,\gamma,\varepsilon,\nu}:=\max\bigg(0\,,\frac{1}{\varepsilon\gamma p_\varepsilon}\(\frac{n}{n-p}\)^k\(1-\frac{\gamma_{i_1,\dotsc,i_k}}{p_*-1+\nu}\)\sum_{i\in\mathcal{I}_1\cup\mathcal{I}_2}\frac{p_\varepsilon-p_i}{p_i}\bigg).
\end{equation}
We separate two cases: 
\renewcommand{\labelitemi}{$-$}
\begin{itemize}
\item Case 1: $p_*-1+\nu\le\gamma_{i_1,\dotsc,i_k}<\frac{n}{p}\(p_\varepsilon-p\)$ (in which case $\tau_{i_1,\dotsc,i_k,\gamma,\varepsilon,\nu}=0$).
\item Case 2: $\gamma_{i_1,\dotsc,i_k}<p_*-1+\nu$ (in which case $\tau_{i_1,\dotsc,i_k,\gamma,\varepsilon,\nu}>0$).
\end{itemize}
We begin with proving \eqref{Lem8Eq8} in Case 1. By interpolation (see, for instance, Grafakos~\cite{Gra}*{Proposition~1.1.14}), and by Lemmas~\ref{Lem1} and~\ref{Lem3}, we obtain
\begin{align}\label{Lem8Eq10}
\left\|u\right\|_{L^{\gamma_{i_1,\dotsc,i_k}}\(A_\varepsilon\(R,\lambda_{k,\gamma,\varepsilon}\)\)}&\le\(\frac{\gamma_{i_1,\dotsc,i_k}}{\gamma_{i_1,\dotsc,i_k}-p_*+1}\)^{\frac{1}{\gamma_{i_1,\dotsc,i_k}}}\left\|u\right\|_{L^{p_*-1,\infty}\(A_\varepsilon\(R,\lambda_{k,\gamma,\varepsilon}\)\)}^{\frac{p_*-1}{\gamma_{i_1,\dotsc,i_k}}}\left\|u\right\|_{L^\infty\(A_\varepsilon\(R,\lambda_{k,\gamma,\varepsilon}\)\)}^{1-\frac{p_*-1}{\gamma_{i_1,\dotsc,i_k}}}\nonumber\\
&\le C\nu^{\frac{-1}{\gamma_{i_1,\dotsc,i_k}}}\le C\nu^{\frac{-1}{p_*-1}}
\end{align}
for some constant $C=C\big(n,\p,K_0,\left\|u\right\|_{L^{p_*-1,\infty}\(\R^n\)}\big)$. Moreover, since $\gamma_{i_1,\dotsc,i_k}<\frac{n}{p}\(p_\varepsilon-p\)$ and $k\le k^+_{\gamma,\varepsilon}$, by \eqref{Lem6Eq2}, we get
\begin{equation}\label{Lem8Eq11}
\frac{\gamma_{i_1,\dotsc,i_k}}{\gamma}\(\frac{n}{n-p}\)^k\le\frac{1}{\varepsilon p_0}\cdot\frac{n}{n-p}\(p_\varepsilon-p\).
\end{equation}
Then \eqref{Lem8Eq8} follows from \eqref{Lem8Eq10} and \eqref{Lem8Eq11}. 

Now, suppose that we are in Case 2. By \eqref{Eq11} and since $\gamma_{i_1,\dotsc,i_{k-1}}\ge\frac{n}{p}\(p_\varepsilon-p\)$ and $p_\varepsilon>p_*$, we obtain
\begin{equation}\label{Lem8Eq12}
\gamma_{i_1,\dotsc,i_k}\ge\frac{n}{p}\(p_\varepsilon-p\)+p_--p_\varepsilon>\frac{n}{p}\(p_*-p\)+p_--p_*=p_--1\,.
\end{equation}
By H\"older's inequality, we then get
\begin{equation}\label{Lem8Eq13}
\left\|u\right\|_{L^{\gamma_{i_1,\dotsc,i_k}}\(A_\varepsilon\(R,\lambda_{k,\gamma,\varepsilon}\)\)}\le\left|A_\varepsilon\(R,\lambda_{k,\gamma,\varepsilon}\)\cap\supp\(u\)\right|^{\frac{1}{\gamma_{i_1,\dotsc,i_k}}-\frac{1}{p_*-1+\nu}}\left\|u\right\|_{L^{p_*-1+\nu}\(A_\varepsilon\(R,\lambda_{k,\gamma,\varepsilon}\)\)}.
\end{equation}
Direct computations yield
\begin{equation}\label{Lem8Eq14}
\left|A_\varepsilon\(R,\lambda_{k,\gamma,\varepsilon}\)\cap\supp\(u\)\right|\le CR^{\frac{1}{\varepsilon p_\varepsilon}\cdot\underset{i\in\mathcal{I}_1\cup\mathcal{I}_2}{\sum}\frac{p_\varepsilon-p_i}{p_i}}
\end{equation}
for some constant $C=C\big(n,\p,R_0\(u\)\big)$, where $R_0\(u\)$ is as in \eqref{Lem8Eq5}. Similarly to \eqref{Lem8Eq10} and \eqref{Lem8Eq11}, we obtain
\begin{equation}\label{Lem8Eq15}
\left\|u\right\|_{L^{p_*-1+\nu}\(A_\varepsilon\(R,\lambda_{k,\gamma,\varepsilon}\)\)}\le C\nu^{\frac{-1}{p_*-1}}
\end{equation}
for some constant $C=C\big(n,\p,K_0,\left\|u\right\|_{L^{p_*-1,\infty}\(\R^n\)}\big)$, and
\begin{equation}\label{Lem8Eq16}
\frac{\gamma_{i_1,\dotsc,i_k}}{\gamma}\(\frac{n}{n-p}\)^k<\frac{1}{\varepsilon p_0}\cdot\frac{p}{n-p}\(p_*-1+\nu\).
\end{equation}
Then \eqref{Lem8Eq8} follows from \eqref{Lem8Eq13}--\eqref{Lem8Eq16}. 

By \eqref{Lem8Eq6} and \eqref{Lem8Eq8}, we obtain
\begin{equation}\label{Lem8Eq17}
\left\|u\right\|_{L^{\gamma}\(A_\varepsilon\(R,\lambda_{0,\gamma,\varepsilon}\)\)}\le\(\widetilde{c}_1c_\nu\)^{\frac{1}{\varepsilon}}\max_{\(i_1,\dotsc,i_k\)\in\Phi_{\gamma,\varepsilon}}R^{\tau_{i_1,\dotsc,i_k,\gamma,\varepsilon,\nu}-\sigma_{i_1,\dotsc,i_k,\gamma,\varepsilon}}
\end{equation}
for all $\nu\in\(0,1\)$, where $\sigma_{i_1,\dotsc,i_k,\gamma,\varepsilon}$ and $\tau_{i_1,\dotsc,i_k,\gamma,\varepsilon,\nu}$ are as in \eqref{Lem8Eq7} and \eqref{Lem8Eq9}. 

We claim that there exists a constant $\nu_0=\nu_0\(n,\p\)$ such that for any $\nu\in\(0,\nu_0\)$, we have
\begin{equation}\label{Lem8Eq18}
\tau_{i_1,\dotsc,i_k,\gamma,\varepsilon,\nu}-\sigma_{i_1,\dotsc,i_k,\gamma,\varepsilon}\le-\frac{1}{\varepsilon p_\varepsilon}\(1-\frac{n}{\gamma p}\(p_\varepsilon-p\)\).
\end{equation}
We prove this claim. By \eqref{Eq12}, we obtain
\begin{align}\label{Lem8Eq19}
\sigma_{i_1,\dotsc,i_k,\gamma,\varepsilon}&=\frac{1}{\varepsilon p_\varepsilon}\bigg(1-\frac{\gamma_{i_1,\dotsc,i_k}}{\gamma}\(\frac{n}{n-p}\)^k+\frac{1}{\gamma}\sum_{j=1}^k\(\frac{n}{n-p}\)^j\(p_\varepsilon-p\)\bigg)\nonumber\\
&=\frac{1}{\varepsilon p_\varepsilon}\bigg(1+\frac{1}{\gamma}\(\frac{n}{n-p}\)^k\(\frac{n}{p}\(p_\varepsilon-p\)-\gamma_{i_1,\dotsc,i_k}\)-\frac{n}{\gamma p}\(p_\varepsilon-p\)\bigg).
\end{align}
In case $p_*-1+\nu\le\gamma_{i_1,\dotsc,i_k}<\frac{n}{p}\(p_\varepsilon-p\)$, since $\tau_{i_1,\dotsc,i_k,\gamma,\varepsilon,\nu}=0$, we deduce \eqref{Lem8Eq18} directly from \eqref{Lem8Eq19}. In the remaining case $\gamma_{i_1,\dotsc,i_k}<p_*-1+\nu$, by \eqref{Lem8Eq9} and \eqref{Lem8Eq19}, we obtain
\begin{multline}\label{Lem8Eq20}
\tau_{i_1,\dotsc,i_k,\gamma,\varepsilon,\nu}-\sigma_{i_1,\dotsc,i_k,\gamma,\varepsilon}\le-\frac{1}{\varepsilon p_\varepsilon}\bigg(1+\frac{1}{\gamma}\(\frac{n}{n-p}\)^k\bigg(\frac{n}{p}\(p_\varepsilon-p\)-\gamma_{i_1,\dotsc,i_k}\\
-\(1-\frac{\gamma_{i_1,\dotsc,i_k}}{p_*-1+\nu}\)\sum_{i\in\mathcal{I}_1\cup\mathcal{I}_2}\frac{p_\varepsilon-p_i}{p_i}\bigg)-\frac{n}{\gamma p}\(p_\varepsilon-p\)\bigg)\,.
\end{multline}
If $\nu$ is small enough so that $p_*-1+\nu<\frac{n}{p}\(p_0-p\)$, i.e. $\nu<\frac{n}{p}\(p_0-p_*\)$, then
\begin{equation}\label{Lem8Eq21}
1-\frac{1}{p_*-1+\nu}\sum_{i\in\mathcal{I}_1\cup\mathcal{I}_2}\frac{p_\varepsilon-p_i}{p_i}<\frac{p}{n\(p_0-p\)}\sum_{i\in\(\mathcal{I}_1\cup\mathcal{I}_2\)^c}\frac{p_0-p_i}{p_i}<0\,,
\end{equation}
where $\(\mathcal{I}_1\cup\mathcal{I}_2\)^c:=\left\{1,\dotsc,n\right\}\backslash\(\mathcal{I}_1\cup\mathcal{I}_2\)$.
It follows from \eqref{Lem8Eq12}, \eqref{Lem8Eq20}, and \eqref{Lem8Eq21} that 
\begin{equation}\label{Lem8Eq22}
\tau_{i_1,\dotsc,i_k,\gamma,\varepsilon,\nu}-\sigma_{i_1,\dotsc,i_k,\gamma,\varepsilon}\le-\frac{1}{\varepsilon p_\varepsilon}\bigg(1+\frac{1}{\gamma}\(\frac{n}{n-p}\)^k\varphi_\nu\(p_\varepsilon\)-\frac{n}{\gamma p}\(p_\varepsilon-p\)\bigg)
\end{equation}
for all $\nu\in(0,\frac{n}{p}\(p_0-p_*\))$, where
\begin{align*}
\varphi_\nu\(q\):&=q-p_--\bigg(1-\frac{\frac{n}{p}\(q-p\)+p_--q}{p_*-1+\nu}\bigg)\underset{i\in\mathcal{I}_1\cup\mathcal{I}_2}{\sum}\frac{q-p_i}{p_i}\\
&=q-p_--\frac{q-p_--\frac{n}{p}\(q-p_*\)}{p_*-1+\nu}\underset{i\in\mathcal{I}_1\cup\mathcal{I}_2}{\sum}\frac{q-p_i}{p_i}
\end{align*}
for all $q\in\R$. By \eqref{Lem8Eq1} and by definition of $\overline{p}_0$, we obtain $\varphi_0\(p_0\)>0$. Moreover, it can easily be seen that $\varphi_0\(p_*\)\le0$. Observing that $\varphi_0$ is a quadratic polynomial with positive leading coefficient, we then get that $\varphi_0$ is increasing in $\[p_0,\infty\)$. By continuity of $\varphi_\nu$ with respect to $\nu$, it follows that $\varphi_\nu\(p_\varepsilon\)\le0$ provided that $\nu<\nu_0$ for some constant $\nu_0=\nu_0\(n,\p\)$. By \eqref{Lem8Eq22}, we then get \eqref{Lem8Eq18}.

Finally, we fix $\nu=\nu_0/2$, and we obtain \eqref{Lem8Eq4} by passing to the limit as $\gamma\to\infty$ into \eqref{Lem8Eq17} and \eqref{Lem8Eq18} and using the fact that $p_\varepsilon>p_0$ and $R>1$. This ends the proof of Lemma~\ref{Lem8}.
\endproof

Now, we can conclude the proof of Theorem~\ref{Th6}.

\proof[Proof of Theorem~\ref{Th6}]
We proceed by contradiction. Suppose that there exists a solution $u$ of \eqref{Eq1} such that 
\begin{equation}\label{Th6Eq1}
p_0:=\max\(\left\{p_i\in\p\,:\,R_i\(u\)=\infty\right\}\)>\overline{p}_0\,,
\end{equation}
where $R_i\(u\)$ is as in \eqref{Lem8Eq2}. Then we can apply Lemma~\ref{Lem8}. For any $\varepsilon\in\(0,1\)$ and $x\in\R^n$, it follows from \eqref{Lem8Eq4} that 
\begin{equation}\label{Th6Eq2}
\left|u\(x\)\right|\le \big(c_2R_\varepsilon\(x\)^{-\frac{1}{p_\varepsilon}}\big)^{\frac{1}{\varepsilon}}\quad\text{where}\quad R_\varepsilon\(x\):=\sum_{i\in\mathcal{I}_2}\left|x_i\right|^{p_\varepsilon}
\end{equation}
provided that $\sum_{i\in\mathcal{I}_1}\left|x_i\right|^{\frac{p_\varepsilon p_i}{p_\varepsilon-p_i}}<\frac{5}{4}R_\varepsilon\(x\)^{1/\varepsilon}$ and $A_\varepsilon\(R_\varepsilon\(x\),1/2\)\cap B_{\p}\(0,\max\(r,1\)\)=\emptyset$, where $\mathcal{I}_1$ and $\mathcal{I}_2$ are as in Lemma~\ref{Lem8}, and $r$ is as in Lemma~\ref{Lem3}. One easily gets that there exists a constant $R_r=R\(n,\p,r\)>1$ such that for any $\varepsilon\in\(0,1\)$ and $R>R_r$, we have $A_\varepsilon\(R,1/2\)\cap B_{\p}\(0,\max\(r,1\)\)=\emptyset$. By passing to the limit as $\varepsilon\to0$ into \eqref{Th6Eq2}, we then obtain that $u\(x\)=0$ for all $x\in\R^n$ such that 
$$\sum_{i\in\mathcal{I}_2}\left|x_i\right|^{p_0}>\max\(R_r,c_2^{\,p_0}\),$$
and hence $R_i\(u\)<\infty$ for all $i\in\mathcal{I}_2$, which is in contradiction with \eqref{Th6Eq1}. This ends the proof of Theorem~\ref{Th6}. 
\endproof

\begin{remark}\label{Rem2}
As one can see from the above proof, the constant $R_0$ that we obtain in Theorem~\ref{Th6} depends on $n$, $\p$, $\Lambda$, $\kappa$, $r$, $r_\kappa\(u\)$, and $\left\|u\right\|_{L^{p_*-1,\infty}\(\R^n\)}$.
\end{remark}

\section{The decay estimates}\label{Sec7}

In this section, we prove Theorem~\ref{Th1} in case $p_+<p_*$ and Theorem~\ref{Th7} below in case $p_*\le p_+<p^*$. The latter implies Theorem~\ref{Th2} in case $p_+=p_*$ and allows us to complete the proof of Theorem~\ref{Th3} in case $p_*<p_+<p^*$. 

\medskip
We let $\overline{p}_0$ and $\mathcal{I}_0$ be as in Section~\ref{Sec6}. We define $q_0$ as the largest real number such that for any $q>q_0$, we have
\begin{equation}\label{Eq17}
\big(q-p_--\frac{n}{p}\(q-p_*\)\big)\sum_{i\in\mathcal{I}_0^c}\frac{q-p_i}{p_i}<\(p_*-1\)\(q-p_-\),
\end{equation}
where $\mathcal{I}_0^c:=\left\{1,\dotsc,n\right\}\backslash\mathcal{I}_0$. It easily follows from the definition of $\overline{p}_0$ and the fact that $\overline{p}_0<p_+$ in case $p_+>p_*$ that
$$\left\{\begin{aligned}
&q_0=\overline{p}_0=p_*&&\text{in case }p_+\le p_*\,,\\
&\overline{p}_0\le q_0<p_+&&\text{in case }p_+>p_*\,.
\end{aligned}\right.$$

\medskip
We prove the following result.

\begin{theorem}\label{Th7}
Assume that $p_*\le p_+<p^ *$. Let $f:\R^n\times\R\to\R$ be a Caratheodory function such that \eqref{Eq2} holds true and $u$ be a solution of \eqref{Eq1}. Let $q_0$ be defined as above. Then for any $q>q_0$, there exists a constant $C_q=C\(n,\p,\Lambda,u,q\)$ such that
\begin{equation}\label{Th7Eq}
\left|u\(x\)\right|^q+\sum_{i=1}^n\left|\partial_{x_i}u\(x\)\right|^{p_i}\le C_q\bigg(1+\sum_{i\in\mathcal{I}_0^c}\left|x_i\right|^{\frac{qp_i}{q-p_i}}\bigg)^{-1}\quad\text{for a.e. }x\in\R^n.
\end{equation}
\end{theorem}

We conclude the proofs of Theorems~\ref{Th2} and~\ref{Th3} as follows.

\proof[Proof of Theorem~\ref{Th2}]
In case $p_+=p_*$, since $q_0=\overline{p}_0=p_*$, we get that \eqref{Th7Eq} holds true for all $q>p_*$. Since in this case we have $\mathcal{I}_0^c=\left\{1,\dotsc,n\right\}$, this is exactly the result in Theorem~\ref{Th2}.
\endproof

\proof[Proof of Theorem~\ref{Th3}]
In case $p_*<p_+<p^ *$, Points (i) and (ii) in Theorem~\ref{Th3} follow directly from Theorems~\ref{Th6} and~\ref{Th7} and the fact that $\overline{p}_0\le q_0<p_+$.
\endproof

Now, it remains to prove Theorems~\ref{Th1} and~\ref{Th7}. By another application of the iteration scheme in Section \ref{Sec5}, we prove the following result.

\begin{lemma}\label{Lem9}
Assume that $p_+<p^*$. Let $f:\R^n\times\R\to\R$ be a Caratheodory function such that \eqref{Eq2} holds true, $u$ be a solution of \eqref{Eq1}, and $\kappa$, $r$, and $K_0$ be as in Lemma~\ref{Lem3}. Let $q=p_*$ in case $p_+<p_*$ and $q\in\(q_0,p^*\)$ in case $p_*\le p_+<p^*$. For any $\lambda\in\(0,1\)$ and $R>1$, we define
\begin{equation}\label{Lem9Eq1}
A_q\(R,\lambda\):=\Omega_{\q}\(\emptyset,1,\mathcal{I}^c_0,R,\lambda\)\quad\text{with}\quad q_i:=\frac{qp_i}{q-p_i}\text{ for all }i\in\mathcal{I}^c_0\,.
\end{equation}
If $A_q\(R,1/2\)\cap B_{\p}\(0,\max\(r,1\)\)=\emptyset$, then
\begin{equation}\label{Lem9Eq2}
\left\|u\right\|_{L^\infty\(A_q\(R,1/4\)\)}\le c_qR^{-\frac{1}{q}},
\end{equation}
for some constant $c_q=c\big(n,\p,\Lambda,K_0,\left\|u\right\|_{L^{p_*-1,\infty}\(\R^n\)},R_0,q\big)$, where $R_0$ is as in Theorem~\ref{Th6}.
\end{lemma}

\proof[Proof of Lemma~\ref{Lem9}]
By Theorem~\ref{Th6}, we obtain that $A_q\(R,\lambda\)\cap\supp\(u\)$ is bounded. By Lemma~\ref{Lem7}, we then get that for any $\varepsilon\in\(0,1\)$, there exists a constant $c_{q,\varepsilon}=c\(n,\p,\Lambda,K_0,q,\varepsilon\)$ such that for any $\gamma>\frac{n}{p}\(\overline{p}_\varepsilon-p\)$, where $\overline{p}_\varepsilon:=\(1+\varepsilon\)\overline{p}_0$, we have
\begin{equation}\label{Lem9Eq3}
\left\|u\right\|_{L^{\gamma}\(A_q\(R,\lambda_{0,\gamma,\varepsilon}\)\)}\le c_{q,\varepsilon}\max_{\(i_1,\dotsc,i_k\)\in\Phi_{\gamma,\varepsilon}}\Big(R^{-\sigma_{i_1,\dotsc,i_k,q,\gamma}}\left\|u\right\|_{L^{\gamma_{i_1,\dotsc,i_k}}\(A_q\(R,\lambda_{k,\gamma,\varepsilon}\)\)}^{\frac{\gamma_{i_1,\dotsc,i_k}}{\gamma}\(\frac{n}{n-p}\)^k}\Big),
\end{equation}
provided that $A_q\(R,1/2\)\cap B_{\p}\(0,\max\(r,1\)\)=\emptyset$, where $\gamma_{i_1,\dotsc,i_k}$ is as in \eqref{Eq11}, $\lambda_{k,\gamma,\varepsilon}$ and $\Phi_{\gamma,\varepsilon}$ are as in \eqref{Lem7Eq2}, and
\begin{equation}\label{Lem9Eq4}
\sigma_{i_1,\dotsc,i_k,q,\gamma}:=\frac{1}{\gamma q}\sum_{j=1}^k\(\frac{n}{n-p}\)^j\(q-p_{i_j}\).
\end{equation}

\proof[End of proof of Lemma~\ref{Lem9} in case $p_*\le p_+<p^*$ and $q_0<q<p^*$]
In this case, we follow in large part the same arguments as in the proof of Lemma~\ref{Lem8}. We set $\varepsilon:=\(q-\overline{p}_0\)/\overline{p}_0$ so that $q=\overline{p}_\varepsilon$. Since $q<p^*$ and $\overline{p}_0\ge p_*$, we get $\varepsilon<\(p^*-p_*\)/p_*\le1$. Similarly to \eqref{Lem8Eq8}, we then obtain that for any $\(i_1,\dotsc,i_k\)\in\Phi_{\gamma,\varepsilon}$ and $\nu\in\(0,1\)$,  there exists a constant $c_\nu=c\big(n,\p,\Lambda,K_0,\left\|u\right\|_{L^{p_*-1,\infty}\(\R^n\)},R_0,\nu\big)$ such that
\begin{equation}\label{Lem9Eq5}
\left\|u\right\|_{L^{\gamma_{i_1,\dotsc,i_k}}\(A_q\(R,\lambda_{k,\gamma,\varepsilon}\)\)}^{\frac{\gamma_{i_1,\dotsc,i_k}}{\gamma}\(\frac{n}{n-p}\)^k}\le c_\nu^{\frac{1}{\varepsilon}}R^{\tau_{i_1,\dotsc,i_k,q,\gamma,\nu}},
\end{equation}
where 
\begin{equation}\label{Lem9Eq6}
\tau_{i_1,\dotsc,i_k,q,\gamma,\nu}:=\max\bigg(0\,,\frac{1}{q\gamma}\(\frac{n}{n-p}\)^k\(1-\frac{\gamma_{i_1,\dotsc,i_k}}{p_*-1+\nu}\)\sum_{i\in\mathcal{I}_0^c}\frac{q-p_i}{p_i}\bigg).
\end{equation}
It follows from \eqref{Lem9Eq3} and \eqref{Lem9Eq5} that 
\begin{equation}\label{Lem9Eq7}
\left\|u\right\|_{L^{\gamma}\(A_q\(R,\lambda_{0,\gamma,\varepsilon}\)\)}\le c_{q,\varepsilon}c_\nu^{\frac{1}{\varepsilon}}\max_{\(i_1,\dotsc,i_k\)\in\Phi_{\gamma,\varepsilon}}R^{\tau_{i_1,\dotsc,i_k,q,\gamma,\nu}-\sigma_{i_1,\dotsc,i_k,q,\gamma}}
\end{equation}
for all $\nu\in\(0,1\)$, where $\sigma_{i_1,\dotsc,i_k,q,\gamma}$ and $\tau_{i_1,\dotsc,i_k,q,\gamma,\nu}$ are as in \eqref{Lem9Eq4} and \eqref{Lem9Eq6}.

In the same way as in the proof of \eqref{Lem8Eq18}, we then obtain 
\begin{equation}\label{Lem9Eq8}
\tau_{i_1,\dotsc,i_k,q,\gamma,\nu}-\sigma_{i_1,\dotsc,i_k,q,\gamma}\le-\frac{1}{q}\(1-\frac{n}{\gamma p}\(q-p\)\)
\end{equation}
provided that
\begin{equation}\label{Lem9Eq9}
q-p_--\frac{q-p_--\frac{n}{p}\(q-p_*\)+\nu}{p_*-1+\nu}\sum_{i\in\mathcal{I}_0^c}\frac{q-p_i}{p_i}>0\,.
\end{equation}
By \eqref{Eq17}, we get that \eqref{Lem9Eq9} holds true provided that $\nu<\nu_0$ for some constant $\nu_0=\nu_0\(n,\p\)$.

Finally, we fix $\nu=\nu_q/2$, and we obtain \eqref{Lem9Eq2} by passing to the limit as $\gamma\to\infty$ into \eqref{Lem9Eq7} and \eqref{Lem9Eq8}. This ends the proof of Lemma~\ref{Lem9}.
\endproof

\proof[Proof of Lemma~\ref{Lem9} in case $p_+<p_*$ and $q=p_*$] In this case, we have $\overline{p}_0=p^*$ and $\mathcal{I}_0^c=\left\{1,\dotsc,n\right\}$. We claim that there exists a constant $\varepsilon_0=\varepsilon_0\(n,\p\)\in\(0,1\)$ such that for any $\varepsilon\in\(0,\varepsilon_0\)$ and $\(i_1,\dotsc,i_k\)\in\Phi_{\gamma,\varepsilon}$, we have
\begin{equation}\label{Lem9Eq10}
\left\|u\right\|_{L^{\gamma_{i_1,\dotsc,i_k}}\(A_q\(R,\lambda_{k,\gamma,\varepsilon}\)\)}^{\frac{\gamma_{i_1,\dotsc,i_k}}{\gamma}\(\frac{n}{n-p}\)^k}\le c_\varepsilon R^{\tau_{i_1,\dotsc,i_k,p_*,\gamma}}
\end{equation}
for some constant $c_\varepsilon=c\big(n,\p,\Lambda,K_0,\left\|u\right\|_{L^{p_*-1,\infty}\(\R^n\)},\varepsilon\big)$, where
\begin{equation}\label{Lem9Eq11}
\tau_{i_1,\dotsc,i_k,p_*,\gamma}:=\frac{p_*-1-\gamma_{i_1,\dotsc,i_k}}{p_*\gamma}\(\frac{n}{n-p}\)^k.
\end{equation}
We assume that $\(1-\varepsilon\)p_*>p$, i.e. $\varepsilon<\frac{p-1}{n-1}$, and we separate two cases: 
\renewcommand{\labelitemi}{$-$}
\begin{itemize}
\item Case 1: $\gamma_{i_1,\dotsc,i_k}\le\frac{n}{p}\(\(1-\varepsilon\)p_*-p\)$,
\item Case 2: $\frac{n}{p}\(\(1-\varepsilon\)p_*-p\)<\gamma_{i_1,\dotsc,i_k}<\frac{n}{p}\(\(1+\varepsilon\)p_*-p\)$.
\end{itemize}
We begin with proving \eqref{Lem9Eq10} in Case 1. By a generalized version of H\"older's inequality (see for instance Grafakos~\cite{Gra}*{Exercise 1.1.11}), we obtain
\begin{multline}\label{Lem9Eq12}
\left\|u\right\|^{\gamma_{i_1,\dotsc,i_k}}_{L^{\gamma_{i_1,\dotsc,i_k}}\(A_q\(R,\lambda_{k,\gamma,\varepsilon}\)\)}\le\frac{p_*-1}{p_*-1-\gamma_{i_1,\dotsc,i_k}}\left|A_q\(R,\lambda_{k,\gamma,\varepsilon}\)\right|^{1-\frac{\gamma_{i_1,\dotsc,i_k}}{p_*-1}}\left\|u\right\|^{\gamma_{i_1,\dotsc,i_k}}_{L^{p_*-1,\infty}\(A_q\(R,\lambda_{k,\gamma,\varepsilon}\)\)}.
\end{multline}
Direct computations give
\begin{equation}\label{Lem9Eq13}
\left|A_q\(R,\lambda_{k,\gamma,\varepsilon}\)\right|\le CR^{\overset{n}{\underset{i=1}{\sum}}\frac{p_*-p_i}{p_*p_i}}=CR^{\frac{p_*-1}{p_*}}
\end{equation}
for some constant $C=C\(n,\p\)$. Since $\gamma_{i_1,\dotsc,i_k}\le\frac{n}{p}\(\(1-\varepsilon\)p_*-p\)$ and $u\in L^{p_*-1,\infty}\(\R^n\)$, it follows from \eqref{Lem9Eq12} and \eqref{Lem9Eq13} that
\begin{equation}\label{Lem9Eq14}
\left\|u\right\|^{\gamma_{i_1,\dotsc,i_k}}_{L^{\gamma_{i_1,\dotsc,i_k}}\(A_q\(R,\lambda_{k,\gamma,\varepsilon}\)\)}\le C\varepsilon^{-1}R^{\frac{p_*-1-\gamma_{i_1,\dotsc,i_k}}{p_*}}
\end{equation}
for some constant $C=C\big(n,\p,\left\|u\right\|_{L^{p_*-1,\infty}\(\R^n\)}\big)$. Moreover, since $k\le k^+_{\gamma,\varepsilon}$, by \eqref{Lem6Eq2}, we get
\begin{equation}\label{Lem9Eq15}
\frac{1}{\gamma}\(\frac{n}{n-p}\)^k<\frac{1}{\varepsilon\(n-1\)}\,.
\end{equation}
Then \eqref{Lem9Eq10} follows from \eqref{Lem9Eq14} and \eqref{Lem9Eq15}. 

Now, suppose that we are in case 2. By interpolation, we obtain
\begin{equation}\label{Lem9Eq16}
\left\|u\right\|_{L^{\gamma_{i_1,\dotsc,i_k}}\(A_q\(R,\lambda_{k,\gamma,\varepsilon}\)\)}\le\left\|u\right\|^\theta_{L^{\frac{n}{p}\(\(1-\varepsilon\)p_*-p\)}\(A_q\(R,\lambda_{k,\gamma,\varepsilon}\)\)}\left\|u\right\|^{1-\theta}_{L^{\frac{n}{p}\(\(1+\varepsilon\)p_*-p\)}\(A_q\(R,\lambda_{k,\gamma,\varepsilon}\)\)},
\end{equation}
where $\theta\in\(0,1\)$ is such that
\begin{equation}\label{Lem9Eq17}
\frac{\theta}{\frac{n}{p}\(\(1-\varepsilon\)p_*-p\)}+\frac{1-\theta}{\frac{n}{p}\(\(1+\varepsilon\)p_*-p\)}=\frac{1}{\gamma_{i_1,\dotsc,i_k}}\,.
\end{equation}
Similarly to \eqref{Lem9Eq14}, we get
\begin{equation}\label{Lem9Eq18}
\left\|u\right\|^{\frac{n}{p}\(\(1-\varepsilon\)p_*-p\)}_{L^{\frac{n}{p}\(\(1-\varepsilon\)p_*-p\)}\(A_q\(R,\lambda_{k,\gamma,\varepsilon}\)\)}\le C\varepsilon^{-1}R^{\frac{n}{p}\varepsilon}
\end{equation}
for some constant $C=C\big(n,\p,K_0,\left\|u\right\|_{L^{p_*-1,\infty}\(\R^n\)}\big)$. On the other hand, Lemma~\ref{Lem4} gives
\begin{equation}\label{Lem9Eq19}
\left\|u\right\|^{\frac{n}{p}\(\(1+\varepsilon\)p_*-p\)}_{L^{\frac{n}{p}\(\(1+\varepsilon\)p_*-p\)}\(A_q\(R,\lambda_{k,\gamma,\varepsilon}\)\)}\le C\max_{i=1,\dotsc,n}\Big(\varepsilon^{-p_i}R^{\frac{p_i-p_*}{p_*}}\left\|u\right\|_{L^{\varepsilon\(n-1\)+p_i-1}\(A_q\(R,1/2\)\)}^{\varepsilon\(n-1\)+p_i-1}\Big)^{\frac{n}{n-p}}
\end{equation}
for some constant $C=C\(n,\p,\Lambda,K_0\)$. We define $\varepsilon_0:=\(p_*-p_+\)/\(p_*+2n-2\)$ so that for any $\varepsilon\in\(0,\varepsilon_0\)$ and $i=1,\dotsc,n$, we have 
$$\varepsilon\(n-1\)+p_i-1<\frac{n}{p}\(\(1-\varepsilon\)p_*-p\).$$
Similarly to \eqref{Lem9Eq14}, we then get
\begin{equation}\label{Lem9Eq20}
\left\|u\right\|^{\varepsilon\(n-1\)+p_i-1}_{L^{\varepsilon\(n-1\)+p_i-1}\(A_q\(R,\lambda_{k+1,\gamma,\varepsilon}\)\)}\le CR^{\frac{p_*-p_i-\varepsilon\(n-1\)}{p_*}}
\end{equation}
for some constant $C=C\big(n,\p,K_0,\left\|u\right\|_{L^{p_*-1,\infty}\(\R^n\)}\big)$. By putting together \eqref{Lem9Eq16}--\eqref{Lem9Eq20}, we obtain
\begin{equation}\label{Lem9Eq21}
\left\|u\right\|^{\gamma_{i_1,\dotsc,i_k}}_{L^{\gamma_{i_1,\dotsc,i_k}}\(A_q\(R,\lambda_{k,\gamma,\varepsilon}\)\)}\le C\varepsilon^{-s}R^{\frac{p_*-1-\gamma_{i_1,\dotsc,i_k}}{p_*}}
\end{equation}
for some constants $C=C\big(n,\p,\Lambda,K_0,\left\|u\right\|_{L^{p_*-1,\infty}\(\R^n\)}\big)$ and $s=s\(n,\p\)>0$. Then \eqref{Lem9Eq10} follows from \eqref{Lem9Eq15} and \eqref{Lem9Eq21}. 

By \eqref{Lem9Eq3} and \eqref{Lem9Eq10}, we obtain that for any $\varepsilon\in\(0,\varepsilon_0\)$, there exists a constant $\widetilde{c}_\varepsilon=c\big(n,\p,\Lambda,K_0,\left\|u\right\|_{L^{p_*-1,\infty}\(\R^n\)},\varepsilon\big)$ such that
\begin{equation}\label{Lem9Eq22}
\left\|u\right\|_{L^\gamma\(A_q\(R,\lambda_{0,\gamma,\varepsilon}\)\)}\le \widetilde{c}_\varepsilon\max_{\(i_1,\dotsc,i_k\)\in\Phi_{\gamma,\varepsilon}}R^{\tau_{i_1,\dotsc,i_k,p_*,\gamma}-\sigma_{i_1,\dotsc,i_k,p_*,\gamma}},
\end{equation}
where $\sigma_{i_1,\dotsc,i_k,p_*,\gamma}$ and $\tau_{i_1,\dotsc,i_k,p_*,\gamma}$ are as in \eqref{Lem9Eq4} and \eqref{Lem9Eq11}.

From \eqref{Eq12}, we derive
$$\sigma_{i_1,\dotsc,i_k,p_*,\gamma}=\frac{1}{p_*}\bigg(1+\frac{1}{\gamma}\(\frac{n}{n-p}\)^k\(p_*-1-\gamma_{i_1,\dotsc,i_k}\)-\frac{p_*-1}{\gamma}\bigg),$$
and hence
\begin{equation}\label{Lem9Eq23}
\tau_{i_1,\dotsc,i_k,p_*,\gamma}-\sigma_{i_1,\dotsc,i_k,p_*,\gamma}=-\frac{1}{p_*}\bigg(1-\frac{p_*-1}{\gamma}\bigg).
\end{equation}
Finally, we fix $\varepsilon=\varepsilon_0/2$, and we obtain \eqref{Lem9Eq2} by passing to the limit as $\gamma\to\infty$ into \eqref{Lem9Eq22} and \eqref{Lem9Eq23}. This ends the proof of Lemma~\ref{Lem9} in case $p_+<p_*$ and $q=p_*$.
\endproof 

Now, we can prove Theorems~\ref{Th7} and~\ref{Th1}.

\proof[Proof of Theorem~\ref{Th7}] As is easily seen, it is sufficient to prove \eqref{Th7Eq} for $q\in\(q_0,p^*\)$. Let $u$ be a solution of \eqref{Eq1} and $q>q_0$. We define 
$$u_R\(y\):=R^{\frac{1}{q}}\cdot u\(\tau_R\(y\)\),\quad\text{where}\quad\tau_R\(y\):=\big(R^{\frac{q-p_1}{qp_1}}y_1,\dotsc,R^{\frac{q-p_n}{qp_n}}y_n\big)$$
for all $R>1$ and $y\in\R^n$. By Lemma~\ref{Lem9}, we obtain
\begin{equation}\label{Th1Eq13}
\left\|u_R\right\|_{L^\infty\(A_q\(1,1/4\)\)}\le c_q\,.
\end{equation}
provided that $A_q\(R,1/2\)\cap B_{\p}\(0,\max\(r,1\)\)=\emptyset$, where $r$ be as in Lemma~\ref{Lem3}. One easily gets the existence of a constant $R_r=R\(n,\p,r\)>1$ such that $A_q\(R,1/2\)\cap B_{\p}\(0,\max\(r,1\)\)=\emptyset$ for all $R>R_r$. Moreover, by \eqref{Eq1}, we obtain
\begin{equation}\label{Th1Eq14}
-\Delta_{\p}u_R=R^{\frac{q-1}{q}}\cdot f\big(\tau_R\(y\),R^{-\frac{1}{q}}\cdot u_R\big)\quad\text{in }\R^n,
\end{equation}
and \eqref{Eq2} gives
\begin{equation}\label{Th1Eq15}
\big|R^{\frac{q-1}{q}}\cdot f\big(\tau_R\(y\),R^{-\frac{1}{q}}\cdot u_R\big)\big|\le\Lambda\cdot R^{\frac{q-p^*}{q}}\cdot\left|u_R\right|^{p^*-1}.
\end{equation}
Since $q-p^*\le0$, by \eqref{Th1Eq13}--\eqref{Th1Eq15} and Lieberman's gradient estimates~\cite{Lie}, we get that there exists a constant $c'_q=c\(n,\p,\Lambda,c_q\)$ such that 
\begin{equation}\label{Th1Eq16}
\left\|\nabla u_R\right\|_{L^\infty\(A_q\(1,1/8\)\)}\le c'_q\,.
\end{equation}
For any $x\in\R^n$, it follows from \eqref{Th1Eq13} and \eqref{Th1Eq16} that
$$\left|u\(x\)\right|^q+\sum_{i=1}^n\left|\partial_{x_i}u\(x\)\right|^{p_i}\le c''_qR\(x\)^{-1}\quad\text{where}\quad R\(x\):=\sum_{i\in\mathcal{I}_0^c}\left|x_i\right|^{\frac{qp_i}{q-p_i}}$$
for some constant $c''_q=c\(n,\p,\Lambda,c_q\)$, provided that $R\(x\)>R_r$. This ends the proof of Theorem~\ref{Th7}.
\endproof

\proof[Proof of Theorem~\ref{Th1}]
We fix $q=p^*$ in this case and we follow the same arguments as in the above proof of Theorem~\ref{Th7}.
\endproof

\begin{remark}\label{Rem3}
As one can see from the above proofs, the constants $C_0$ and $C_q$ that we obtain in \eqref{Th1Eq1} and~\eqref{Th7Eq} depend on $n$, $\p$, $q$, $\Lambda$, $\kappa$, $r$, $r_\kappa\(u\)$, $\left\|u\right\|_{L^{p_*-1,\infty}\(\R^n\)}$, and $\left\|u\right\|_{W^{1,\infty}\(\R^n\backslash\Omega_r\)}$, where $\Omega_r:=\big\{x\in\R^n:\,\sum_{i=1}^n\left|x_i\right|^{\frac{qp_i}{q-p_i}}>R_r\big\}$ for some constant $R_r=R\(n,\p,r\)$.
\end{remark}

\appendix

\section{Kato-type inequality}\label{App}

In this section, we prove a weak version of Kato's inequality~\cite{Kato} for the operator $\Delta_{\p}$. This result is used in Sections~\ref{Sec3} and~\ref{Sec4}. A similar result has been proven by Cuesta~Leon~\cite{Cue} in the context of the $p$--Laplace operator. 

\medskip
For any $f\in L^1_{\loc}\(\R^n\)$, we say that a function $u\in D^{1,\p}\(\R^n\)$ is a solution of the inequality
$$-\Delta_{\p}u\le f\quad\text{in }\R^n$$
if we have
$$\sum_{i=1}^n\int_{\R^n}\left|\partial_{x_i}u\right|^{p_i-2}\(\partial_{x_i}u\)\(\partial_{x_i}\varphi\)dx\le\int_{\R^n}f\varphi\,dx$$
for all nonnegative, smooth function $\varphi$ with compact support in $\R^n$.

\medskip
We state our result as follows.

\begin{proposition}\label{Pr2}
Let $f_1,f_2\in L^1_{\loc}\(\R^n\)$ and $u_1,u_2\in D^{1,\p}\(\R^n\)$ be solutions of the inequalities
\begin{equation}\label{Pr2Eq1}
-\Delta_{\p}u_j\le f_j\quad\text{in }\R^n
\end{equation}
for $j=1,2$. Then the function $u:=\max\(u_1,u_2\)$ is a solution of the inequality
\begin{equation}\label{Pr2Eq2}
-\Delta_{\p}u\le f\quad\text{in }\R^n,
\end{equation}
where $f\(x\):=f_1\(x\)$ if $u_1\(x\)>u_2\(x\)$, $f\(x\):=f_2\(x\)$ if $u_1\(x\)\le u_2\(x\)$ for all $x\in\R^n$. 
\end{proposition}

\proof[Proof of Proposition~\ref{Pr2}]
We essentially follow the lines of Cuesta~Leon~\cite{Cue}*{Proposition~3.2}. For any $\varepsilon>0$ and $x\in\R^n$, we define
$$\eta_{1,\varepsilon}\(x\):=\eta_\varepsilon\(u_1\(x\)-u_2\(x\)\)\quad\text{and}\quad\eta_{2,\varepsilon}\(x\):=1-\eta_{1,\varepsilon}\(x\),$$
where $\eta_\varepsilon\in C^1\(R\)$ is such that $\eta_\varepsilon\equiv0$ in $\(-\infty,0\]$, $\eta_\varepsilon\equiv1$ in $\[1,\infty\)$, $0\le\eta_\varepsilon\le1$ and $\eta'_\varepsilon\ge0$ in $\(0,1\)$. In particular, for $j=1,2$, we have $\eta_{j,\varepsilon}\in D^{1,\p}\(\R^n\)$, $0\le\eta_{j,\varepsilon}\le1$ in $\R^n$, and
\begin{equation}\label{Pr2Eq3}
\eta_{j,\varepsilon}\(x\)\longrightarrow\left\{\begin{aligned}&1&&\text{if }x\in\Omega_j\\&0&&\text{if }x\in\R^n\backslash\Omega_j\end{aligned}\right.
\end{equation}
as $\varepsilon\to0$ for all $x\in\R^n$, where $\Omega_1:=\left\{x\in\R^n:\,u_1\(x\)>u_2\(x\)\right\}$ and $\Omega_2:=\R^n\backslash\Omega_1$. For any nonnegative, smooth function $\varphi$ with compact support in $\R^n$, testing \eqref{Pr2Eq1} with $\varphi\eta_{j,\varepsilon}$ gives
\begin{multline}\label{Pr2Eq4}
\sum_{i=1}^n\int_{\R^n}\left|\partial_{x_i}u_j\right|^ {p_i-2}\(\partial_{x_i}u_j\)\(\partial_{x_i}\varphi\)\eta_{j,\varepsilon}dx\\
+\(-1\)^{j-1}\sum_{i=1}^n\int_{\R^n}\left|\partial_{x_i}u_j\right|^ {p_i-2}\(\partial_{x_i}u_j\)\(\partial_{x_i}u_1-\partial_{x_i}u_2\)\eta'_\varepsilon\(u_1-u_2\)\varphi\,dx\le\int_{\R^n}f_j\varphi\eta_{j,\varepsilon}dx\,.
\end{multline}
By \eqref{Pr2Eq3} and since $f_j\in L^1_{\loc}\(\R^n\)$ and $u_j\in D^{1,\p}\(\R^n\)$, we obtain
\begin{align}
\int_{\R^n}f_j\varphi\eta_{j,\varepsilon}dx&\longrightarrow\int_{\Omega_j}f_j\varphi\,dx\,,\label{Pr2Eq5}\\
\int_{\R^n}\left|\partial_{x_i}u_j\right|^ {p_i-2}\(\partial_{x_i}u_j\)\(\partial_{x_i}\varphi\)\eta_{j,\varepsilon}dx&\longrightarrow\int_{\Omega_j}\left|\partial_{x_i}u_j\right|^ {p_i-2}\(\partial_{x_i}u_j\)\(\partial_{x_i}\varphi\)dx\label{Pr2Eq6}
\end{align}
as $\varepsilon\to0$ for all $i=1,\dotsc,n$. Moreover, since $\eta'_\varepsilon\ge0$ and $\varphi\ge0$, we get
\begin{equation}\label{Pr2Eq7}
\sum_{i=1}^n\int_{\R^n}\(\left|\partial_{x_i}u_1\right|^ {p_i-2}\partial_{x_i}u_1-\left|\partial_{x_i}u_2\right|^ {p_i-2}\partial_{x_i}u_2\)\(\partial_{x_i}u_1-\partial_{x_i}u_2\)\eta'_\varepsilon\(u_1-u_2\)\varphi\,dx\ge0\,.
\end{equation}
It follows from \eqref{Pr2Eq4}--\eqref{Pr2Eq7} that
\begin{equation}\label{Pr2Eq8}
\sum_{j=1}^2\sum_{i=1}^n\int_{\Omega_j}\left|\partial_{x_i}u_j\right|^ {p_i-2}\(\partial_{x_i}u_j\)\(\partial_{x_i}\varphi\)dx\le\sum_{j=1}^2\int_{\Omega_j}f_j\varphi\,dx
\end{equation}
and hence \eqref{Pr2Eq2} holds true since $u\equiv u_j$ and $f\equiv f_j$ on $\Omega_j$ for $j=1,2$.
\endproof

\end{document}